\documentclass[11pt]{article}
\usepackage{latexsym,graphicx}

\font\tenmath=msbm10 scaled 1200
\font\sevenmath=msbm7 scaled 1200
\font\fivemath=msbm5 scaled 1200

\newfam\mathfam \textfont\mathfam=\tenmath
\scriptfont\mathfam=\sevenmath \scriptscriptfont\mathfam=\fivemath
\def\math{\fam\mathfam}
\def\R{{\math R}}
\def\N{{\math N}}
\def\E{{\math E}}

\def\P{{\math P}}

\def \^#1{\if#1i{\accent"5E\i}\else{\accent"5E#1}\fi}

\def \ind {\mbox{\bf 1}}

\def \vfi{\varphi}

\def \e{\varepsilon}
\def \g{\gamma}

\def \D{\Delta}

\def \cqfd{\quad_\diamondsuit} 

\def \ss{\smallskip}
\def \bs{\bigskip}
\def \ni{\noindent}
\def \F{{\cal F}}

\newtheorem{Thm}{Theorem}
\newtheorem{Lem}{Lemma}
\newtheorem{Pro}{Proposition}

\newtheorem{rem}{Remark}
\newenvironment{Rem}{\begin{rem}\rm}{\end{rem}}
\author{\sc Damien Lamberton \thanks{Laboratoire d'analyse et de math\'ematiques appliqu\'ees,
UMR~8050, Univ. Marne-la-Vall\'ee, Cit\'e
Descartes, 5, Bld Descartes, Champs-sur-Marne, F-77454 Marne-la-Vall\'ee Cedex 2, France.  {\tt
damien.lamberton@univ-mlv.fr}} 
\quad Gilles Pag\`es
\thanks{Laboratoire de probabilit\'es et mod\`eles al\'eatoires, UMR~7599, Univ. Paris 6, case 188, 4,
pl. Jussieu, F-75252 Paris Cedex 5. {\tt  gpa@ccr.jussieu.fr}} } 

\title{\bf A penalized bandit algorithm~\thanks{This work has benefitted from the stay of both authors
at the Isaac Newton Institute (Cambridge University) on the program {\em Developments in Quantitative Finance}.}} 

\date{} 

\begin{document}

\maketitle
\begin{abstract}
We study a two armed-bandit algorithm with penalty. We show the 
convergence of the algorithm and establish the rate of convergence.
For some choices of the parameters, we obtain a central limit 
theorem in which the limit distribution is characterized as
the unique stationary distribution of a discontinuous Markov process.
\end{abstract}

\bigskip


\bs \ni {\em 2001 AMS classification:} 62L20,
%Stochastic Approximation
secondary 93C40,
%Adaptive control (code de Narendra lui-meme) 
91E40, 
% memory and learning
68T05, 
%learning and adaptive systems
91B32
%resource and cost allocation
.

 \section*{Introduction} 
In a recent joint work with P. Tarr\`es (see \cite{LAPATA}), we 
studied the convergence  of the so-called two armed bandit 
algorithm. The purpose of the present paper is to investigate a 
modified version of this algorithm, in which a  penalization
is introduced. In the terminology of learning theory
(see~\cite{NA1,NA2}), the algorithm studied in~\cite{LAPATA}
was a Linear Reward-Inaction ($LRI$) scheme, whereas the one we want to 
introduce is a Linear Reward-Penalty ($LRP$) procedure.

In our previous paper, the algorithm  was introduced in a financial 
context as a procedure for the optimal allocation of a fund
between two traders who manage it. 
Imagine that the owner of a fund can share his wealth between two 
traders, say $A$ and $B$, and that, every day, he can evaluate the 
results of one of the traders and, subsequently,
modify 
the percentage of the fund managed by both traders. Denote by 
$X_n$ the percentage managed by trader $A$ at time $n$ ($X_n\!\in [0,1]$). We assume
that the owner selects the trader to be evaluated at random,
in such a way that the probability that $A$ is evaluated at time
$n$ is $X_n$, in order to select preferably the trader in charge of 
the greater part of the fund. In the LRI scheme, if the evaluated 
trader performs well, its share is increased by a fraction
$\gamma_n \!\in(0,1)$  of the share of the other trader,  and nothing happens if
the evaluated trader performs badly. Therefore, the dynamics of the sequence
$(X_n)_{n\geq 0}$ can be modelled as follows:
\[
X_{n+1}=X_n+\gamma_{n+1}\left(
 \ind_{\{U_{n+1}\leq X_n\}\cap A_{n+1}}(1-X_n)
          -\ind_{\{U_{n+1}> X_n\}\cap B_{n+1}}X_n\right),\quad X_0=x\!\in[0,1],
\]
where $(U_n)_{n\geq 1}$ is an iid sequence 
of uniform random variables on the interval $[0,1]$,
$A_n$ (resp. $B_n$) is the event ``trader $A$ (resp. trader $B$) performs well at 
time $n$". We assume $\P(A_n)=p_{\!_A}$, $\P(B_n)=p_{\!_B}$, for $n\geq 1$, 
with $p_{\!_A},p_{\!_B}\in(0,1)$, and independence between these events and the 
sequence $(U_n)_{n\geq 1}$. The point is that the owner of the fund does not know the parameters
$p_{\!_A},\,p_{\!_B}$.

This recursive learning procedure   has been designed in order to assign
asymptotically the whole fund to the best trader. This  means that, if say $p_{\!_A}> p_{\!_B}$,   $X_n$
converges to $1$ with probability $1$ provided $X_0\!\in(0,1)$ (if $p_{\!_A}< p_{\!_B}$, the limit is $0$ with
symmetric results). However this ``infallibility" property needs  some very stringent assumptions on the reward
parameter
$\g_n$ (see~\cite{LAPATA}). Furthermore, the rate of convergence of the procedure either toward its ``target"
$1$ or its ``trap" $0$ is not ruled by a
CLT with rate $\sqrt{\g_n}$ like standard stochastic approximation algorithms (see~\cite{KUYI}). It is shown
in~\cite{LAPA} that this rate is quite non-standard, strongly depends on the  (unknown)  values $p_{\!_A}$ and
$p_{\!_B}$ and becomes very poor as  these probabilities get close to each other.

\smallskip
In order to improve the efficiency of the algorithm,
one may imagine to introduce a penalty when an evaluated  trader
 has unsatisfactory performances. More precisely, if
 the evaluated trader at time $n$  performs badly, its share is decreased by
 a    penalty  factor $\rho_n\gamma_n$. This leads   to the following $LRP$ 
 --~or ``penalized two-armed bandit~--  procedure
\begin{eqnarray*}
 X_{n+1}&=&X_n+\gamma_{n+1}\left(
 \ind_{\{U_{n+1}\leq X_n\}\cap A_{n+1}}(1-X_n)
          -\ind_{\{U_{n+1}> X_n\}\cap B_{n+1}}X_n\right)\\
     &&	-\gamma_{n+1}\rho_{n+1}\left(
	X_n\ind_{\{U_{n+1}\leq X_n\}\cap A^c_{n+1}}
	-(1-X_n)\ind_{\{U_{n+1}>X_n\}\cap B^c_{n+1}}
	\right), \quad n\in \N,
\end{eqnarray*}
where the notation $A^c$ is used for the complement of an event $A$.
The precise assumptions on the reward rate $\gamma_n$ and the penalty
rate   $\g_n\rho_n$ will be given in the following sections.

The paper is organized as follows. In  Section~\ref{Convergence},
we discuss the convergence of the sequence $(X_n)_{n\geq 0}$. First  we show that, if  $\rho_n$ is a positive 
constant $\rho$,
the sequence converges with probability one to a limit
$x^*_\rho\in(0,1)$  
 satisfying $x^*_\rho>\frac 12$ if and only if $p_{\!_A}>p_{\!_B}$, so that, although the algorithm manages to 
distinguish which trader is better, it does not assign
the whole fund to the best trader.  To get rid of this limitation,  we consider
a sequence $(\rho_n)_{n\geq 1}$ which goes to zero  so  that the penalty
rate becomes negligible with respect to the reward rate ($\g_n\rho_n
=o(\g_n)$). This framework seems  new  in the learning theory literature.
Then, we are able to show  that  the algorithm is infallible $i.e.$, if
$p_{\!_A}>p_{\!_B}$, then 
$\displaystyle \lim_{n\to\infty}X_n=1$ almost surely, under very  light  conditions on
the reward rate $\gamma_n$ (and $\rho_n$). From a stochastic approximation
viewpoint, this modification of the original procedure has the same mean
function and time scale (hence the same target and trap,
see~(\ref{canonic})) but it  always keeps the algorithm away from the
trap without adding noise at these equilibria. In fact, it was necessary
not to add noise at these points in order to remain inside the domain
$[0,1]$. 

The other two sections are devoted to the rate of convergence.
In Section~\ref{SectionPtw}, we show that under some conditions 
(including $\displaystyle \lim_{n\to\infty}\gamma_n/\rho_n=0$)
the sequence $Y_n = (1-X_n)/\rho_n$ converges in probability to 
$(1-p_{\!_A})/\pi$, where $\pi=p_{\!_A}-p_{\!_B}>0$. With additional assumptions,
we prove that this convergence   occurs with probability $1$.
In Section~\ref{Weak}, we show that if the ratio $\gamma_n/\rho_n$
goes to a positive limit as $n$ goes to infinity, then $(Y_n)_{n\ge 1}$ converges  in a weak
sense  to  a probability distribution $\nu$. This
distribution is identified as the unique  stationary distribution of a discontinuous Markov process. This result
is obtained by using weak functional methods applied to a  re-scaling of the algorithm. This approach  can be
seen as an extension of the $SDE$ method used to prove the CLT in a more standard framework of stochastic
approximation (see~\cite{KUYI}). Furthermore, we show that $\nu$ is absolutely continuous with continuous,
possibly non-smooth,  piecewise ${\cal C}^\infty$ density.  An interesting  consequence of these results for practical
applications is that, by choosing
$\rho_n$ and
$\gamma_n$ proportional to
$n^{-1/2}$, one can achieve {\em convergence at the rate $1/\sqrt{n}$, without
any    a priori   knowledge about the values of $p_{\!_A}$ and $p_{\!_B}$}. This is in contrast
with the case of the LRI procedure, where the rate of convergence 
depends heavily on these parameters (see~\cite{LAPA}) and becomes quite poor when they get close to
each other.

\bs
\noindent {\sc Notation.}  
 Let $(a_n)_{n\ge 0}$ and $(b_n)_{n\ge 0}$ be two sequences of positive real numbers. The symbol
$a_n\sim b_n$ means $ a_n=b_n +o(b_n)$.

 \section{Convergence of the LRP algorithm}
 \label{Convergence}
 \subsection{Some classical background on stochastic
 approximation} We will rely on the
 $ODE$ lemma recalled below for a stochastic procedure
 $(Z_n)$ taking its values in a given compact interval
 $I$.
 \begin{Thm} \label{ThmODE} $(a)$ {\sc Kushner \& Clark's $ODE$ Lemma~(see~\cite{KUCL}):}  Let $g:I\to \R$ such that $Id+g$ leaves $I$ stable
 (\footnote{then for every
 $\g\!\in[0,1]$, $Id+\g g= \g(Id+g)+(1-\g)Id $ still takes values in the convex set $I$}). Then, consider the recursively defined stochastic
 approximation procedure  defined  on $I$ by
 \[
 Z_{n+1}=Z_n+\g_{n+1}(g(Z_n) + \Delta R_{n+1}),\;n\ge 0, \quad Z_0\!\in I,
 \]
 where $(\g_n)_{n\ge 1}$ is a sequence of $[0,1]$-valued real numbers satisfying $\g_n\to 0$ and $\sum_{n\ge 1}\g_n =+\infty$.  Set $N(t):=
 \min\{n\,:\,
 \g_1+\cdots+\g_{n+1} > t\}$. If, for every
 $T>0$,
 \begin{equation}\label{ODEcond}
 \max_{N(t) \le n\le N(t+T)}\left|\sum_{k=N(t)+1}^n\g_k\,\Delta R_k\right|\longrightarrow  0\qquad \P\mbox{-a.s. as} \quad t\to
 +\infty.
 \end{equation}
 Let $z^*$ be an attracting zero of $g$ in $I$ and $G(z^*)$ its attracting interval. Then, on the event
 \[
 \{Z_n \mbox{ visits infinitely often a compact subset of } G(z^*)\}\qquad Z_n \mathop{\longrightarrow}^{a.s.}z^*.
 \]

 \noindent $(b)$ {\sc The Hoeffding condition~(see~\cite{BE}):} If $(\Delta R_n)_{n\ge 0}$ is a sequence of $L^\infty$-bounded martingale
 increments, if $(\g_n)$ is nonincreasing and $\displaystyle\sum_{n\ge 1} e^{-\frac{\vartheta}{\g_n}}<+\infty$ for every $\vartheta>0$, then
 Assumption~(\ref{ODEcond}) is satisfied. \end{Thm}

 \bigskip
 \noindent {\bf Remark.} The monotonous assumption on the sequence $\g$ can be relaxed into $\g_n\to 0$ and $\displaystyle
 \sup_{n,k\ge 1}\frac{\g_{n+k}}{\g_n}<+\infty$
 \subsection{Basic properties of the LRP algorithm}
 We first recall the definition of the algorithm. We are interested 
 in the asymptotic behavior of the sequence $(X_n)_{n\in\N}$, where 
 $X_0=x$, with $x\in(0,1)$, and
 \begin{eqnarray*}
 X_{n+1}&=&X_n+\gamma_{n+1}\left(
 \ind_{\{U_{n+1}\leq X_n\}\cap A_{n+1}}(1-X_n)
          -\ind_{\{U_{n+1}> X_n\}\cap B_{n+1}}X_n\right)\\
     &&	-\gamma_{n+1}\rho_{n+1}\left(
	X_n\ind_{\{U_{n+1}\leq X_n\}\cap A^c_{n+1}}
	-(1-X_n)\ind_{\{U_{n+1}>X_n\}\cap B^c_{n+1}}
	\right), \quad n\in \N.
\end{eqnarray*}
Throughout the paper, we assume that
$(\gamma_n)_{n\geq 1}$ is a non-increasing sequence of positive numbers 
satisfying $\gamma_n<1$, $\displaystyle \sum_{n=1}^\infty \gamma_n=+\infty$ and
\[
\forall \vartheta>0,\quad \sum_ne^{-\frac{\vartheta}{\gamma_n}}<\infty, 
\]
and that $(\rho_n)_{n\geq 1}$ is a sequence of positive numbers satisfying 
$\gamma_n\rho_n<1$; $(U_n)_{n\geq 1}$ is 
a sequence of independent random variables which are uniformly 
distributed on the interval $[0,1]$, the events $A_n$, $B_n$
satisfy
\[
\P(A_n)=p_{\!_A},\quad \P(B_n)=p_{\!_B}, \quad n\in\N,
\]
where $0< p_{\!_B}\leq p_{\!_A}<1$, and the sequences $(U_n)_{n\geq 1}$
and $(\ind_{A_n}, \ind_{B_n})_{n\geq 1}$ are independent.
The natural filtration of the sequence $(U_n,\ind_{A_n}, \ind_{B_n})_{n\geq 1}$
is denoted by $(\F_n)_{n\geq 0}$ and
we set
\[
\pi=p_{\!_A}-p_{\!_B}.
\]
With this notation, we have, for $n\geq 0$,
\begin{equation}\label{mean-algo}
X_{n+1}=X_n+\gamma_{n+1}\left(\pi h(X_n)+\rho_{n+1}\kappa(X_n)\right)+\gamma_{n+1}\Delta M_{n+1}, 
\end{equation}
where the functions $h$ and $\kappa$ are defined by
\[
h(x)=x(1-x), \quad \kappa(x)=-(1-p_{\!_A})x^2+(1-p_{\!_B})(1-x)^2, \quad 0\leq x\leq 1,
\]
$\Delta M_{n+1}=M_{n+1}-M_n$, and the sequence $(M_n)_{n\geq 0}$ 
is the martingale defined by $M_0=0$ and
\begin{eqnarray}
\Delta M_{n+1}&=&\ind_{\{U_{n+1}\leq X_n\}\cap A_{n+1}}(1-X_n)
          -\ind_{\{U_{n+1}> X_n\}\cap B_{n+1}}X_n-\pi 
          h(X_n)\nonumber\\
&&  -\rho_{n+1}\left(
	X_n\ind_{\{U_{n+1}\leq X_n\}\cap A^c_{n+1}}
	-(1-X_n)\ind_{\{U_{n+1}>X_n\}\cap B^c_{n+1}}+\kappa(X_n)\right).
	\label{M-dynamics}
\end{eqnarray}
Observe that the increments $\Delta M_{n+1}$ are bounded.
\subsection{Constant penalty rate}
In this subsection, we assume
\[
\forall n\geq 1, \quad \rho_n=\rho,
\]
with $0<\rho\leq 1$.
We then have
\[
X_{n+1}=X_n+\gamma_{n+1}\left(h_\rho(X_n)+\Delta M_{n+1}\right),
\]
where 
\[
h_\rho(x)=\pi h(x)+\rho\kappa(x),\quad 0\leq x\leq 1.
\]
Note that $h_\rho(0)=\rho(1-p_{\!_B})>0$ and $h_\rho(1)=-\rho(1-p_{\!_A})<0$,
and that there exists a unique $x^*_\rho\in (0,1)$ such that 
$h_\rho(x^*_\rho)=0$. By a straightforward computation,
we have
\begin{eqnarray*}
 x^*_\rho &=& \frac{\pi -2\rho(1-p_{\!_B}) +\sqrt{\pi^2 +4\rho^2(1-p_{B}) (1-p_{A})}}{2\pi (1-\rho)}\quad\mbox{ if }\; \pi \neq 0
 \mbox{ and } \rho\neq 1\\
 *[.5em] &=& \frac{(1-p_{\!_A}) }{(1-p_{\!_A}) +(1-p_{\!_B}) }\qquad\mbox{ if }\; \pi= 0\;\mbox{ or }\; \rho= 1.
 \end{eqnarray*}
 In particular, $x^*_\rho =1/2$ if $\pi=0$ regardless of the value of $\rho$. 
 We also have $h_\rho(1/2)=\pi(1+\rho)/4\geq 0$,
 so that
 \begin{equation}\label{xrhoet1/2}
 x^*_\rho>1/2 \quad\mbox{ if }\quad \pi>0.
 \end{equation}
Now, let $x$ be a solution of the ODE $dx/dt=h_\rho(x)$.
If $x(0)\in [0,x^*_\rho]$, $x$ is non-decreasing and 
$\displaystyle \lim_{t\to\infty}x(t)=x^*_\rho$. If $x(0)\in [x^*_\rho,1]$, $x$ is 
non-increasing and 
$\displaystyle \lim_{t\to\infty}x(t)=x^*_\rho$. It follows that the interval $[0,1]$
is a domain of attraction for $x^*_\rho$.
Consequently, using Kushner and Clark's  $ODE$ Lemma (see Theorem~\ref{ThmODE}), one reaches the following conclusion.
 \begin{Pro} Assume that $\rho_n =\rho \!\in(0,1]$, then
 \[
 X_n \mathop{\longrightarrow}^{a.s.} x^*_\rho \qquad \mbox{ as } n\to \infty.
 \]
 \end{Pro}
 The natural interpretation, given the above inequalities on $x^*_\rho$, is that this algorithm never fails in pointing
 the best trader  thanks to Inequality~(\ref{xrhoet1/2}), but it never 
 assigns the whole fund to this trader as the original $LRI$
 procedure did.
\subsection{Convergence when the penalty rate goes to zero}
\begin{Pro}\label{Lim=0ou1}
    Assume $\displaystyle \lim_{n\to \infty}\rho_n=0$.
    The sequence $(X_n)_{n\in\N}$ is 
    almost surely convergent and its limit $X_\infty$ satisfies
    $X_\infty\in\{0,1\}$ with probability 1.
\end{Pro}
\noindent {\sc Proof:}  We first    write  the algorithm in its canonical  form
 \begin{equation}\label{canonic}
 X_{n+1} = X_n +\g_{n+1}(\pi \, h(X_n) + \Delta R_{n+1})\qquad \mbox{ with } 
 \qquad \Delta R_n = \Delta M_n +\rho_n \kappa(X_{n-1}).
 \end{equation}
  It is straightforward to check that the $ODE$ $\dot{x}= h(x)$ has two equilibrium points,   $0$ and $1$, $1$ being attractive
 with $(0,1]$ as an attracting interval and $0$ is unstable.

  Since the martingale increments $\Delta M_n$ are bounded,  it  follows  
  from the assumptions on the sequence $(\gamma_n)_{n\geq 1}$ and  the Hoeffding condition  (see
 Theorem~\ref{ThmODE}$(b)$)   that
 \[
 \max_{N(t) \le n\le N(t+T)}|\sum_{k=N(t)+1}^n\g_k\Delta M_{k}|\stackrel{\P\mbox{-}a.s.}{\longrightarrow}0\quad \mbox{ as }t\to +\infty
 \]
 for every $T>0$. On the other hand  the  function $\kappa$ being bounded 
 on $[0,1]$ and $\rho_n$ converging to $0$, we have, for every $T>0$,
 \[
 \max_{N(t) \le n\le N(t+T)}|\sum_{k=N(t)+1}^n\g_k\rho_k \kappa(X_{k-1})|\le \|k\|_{[0,1]} (T+\g_{_{N(t+T)}})\max_{k\ge N(t)+1}\!\!\!\!\rho_k\longrightarrow 0\quad \mbox{ as }t\to
 +\infty.
 \]
 Finally,  the sequence $(\Delta R_n)_{n\ge 1}$ satisfies Assumption~(\ref{ODEcond}).
Consequently, either $X_n$ visits infinitely often an interval $[\varepsilon,1]$ for some $\varepsilon>0$ and $X_n $ converges toward
 $1$, or $X_n$ converges toward $0$.
$\cqfd$
\begin{Rem}\label{Rem1}
    If $\pi=0$, $i.e.$ $p_{\!_A}=p_{\!_B}$, the algorithm reduces to
    \[
     X_{n+1}=X_n+\gamma_{n+1}\rho_{n+1}(1-p_{\!_A})(1-2X_n)+
     \gamma_{n+1}\Delta M_{n+1}.
    \]
    The number $1/2$ is the unique equilibrium of the ODE 
    $\dot{x}=(1-p_{\!_A})(1-2x)$, and the interval $[0,1]$ is a domain
    of attraction. Assuming $\sum_{n=1}^\infty 
    \rho_n\gamma_n=+\infty$, and that the sequence 
    $(\gamma_n/\rho_n)_{n\geq 1}$ is non-increasing and satisfies
    \[
    \forall \vartheta>0,\quad \sum_{n=1}^\infty 
    \exp\left(-\vartheta\frac{\rho_n}{\gamma_n}\right)<+\infty,
    \]
    it can be proved, using the Kushner-Clark $ODE$ Lemma
(Theorem~\ref{ThmODE}), that 
    $\displaystyle \lim_{n\to\infty}X_n=1/2$ almost surely. As concerns
the asymptotics of the algorithm when $\pi=0$ and 
$\g_n=g\,\rho_n$ (for which the above condition is not satisfied), we
refer to the final remark of the paper. 
\end{Rem}

From now on, we will assume that $p_{\!_A}>p_{\!_B}$. The next
proposition shows that the penalized algorithm is infallible under very
light assumptions on $\g_n$ and $\rho_n$. 

\begin{Pro}\label{lim=0} (Infallibility)
    Assume $\displaystyle \lim_{n\to \infty}\rho_n=0$.
    If the sequence $(\gamma_n/\rho_n)_{n\geq 1}$ is bounded and 
    $\sum_n\gamma_n\rho_n=\infty$, and if $\pi>0$,
      we have $\displaystyle \lim_{n\to \infty}X_n=1$ almost surely.
    \end{Pro}
\noindent {\sc Proof:} 
We have from (\ref{mean-algo}), since $h\geq 0$ on the interval $[0,1]$,
\[
X_n\geq 
X_0+\sum_{j=1}^n\gamma_j\rho_j\kappa(X_{j-1})+
\sum_{j=1}^n\gamma_j\D M_j, \quad n\geq 1.
\]
Since the jumps $\D M_j$ are bounded, we have
\[
\left|\left|\sum_{j=1}^n\gamma_j\D M_j\right|\right|^2_{L^2}\leq 
C\sum_{j=1}^n\gamma_j^2\leq 
C\sup_{j\in\N}(\gamma_j/\rho_j)\sum_{j=1}^n\gamma_j\rho_j,
\]
for some positive  constant $C$.
Therefore, since $\sum_n\gamma_n\rho_n=\infty$,
\[
L^2\mbox{-}\lim_{n\to\infty}\frac{\sum_{j=1}^n\gamma_j\D 
M_j}{\sum_{j=1}^n\gamma_j\rho_j}=0 \quad \mbox{ so that  }\quad
\limsup_n \frac{\sum_{j=1}^n\gamma_j\D 
M_j}{\sum_{j=1}^n\gamma_j\rho_j}\ge 0 \qquad a.s..
\]
Now, on the set $\{X_\infty=0\}$, we have
\[
\lim_{n\to\infty}
\frac{\displaystyle\sum_{j=1}^n\gamma_j\rho_j\kappa(X_{j-1})}{\displaystyle\sum_{j=1}^n\gamma_j\rho_j}=\kappa(0)>0.
\]
Hence, it follows that, still on the set $\{X_\infty=0\}$,
\[
\limsup_{n\to\infty}\frac{X_n}{\displaystyle\sum_{j=1}^n\gamma_j\rho_j}>0.
\]
Therefore, we must have  $\P(X_\infty=0)=0$.$\cqfd$

\bs

The following Proposition will give a control on  
the conditional variance process of the martingale $(M_n)_{n\in\N}$ which
will be crucial to elucidate the rate of convergence of the algorithm.
\begin{Pro}\label{ProVarCond}
    We have, for $n\geq 0$,
    \[
     \E\left(\D M_{n+1}^2\;|\;\F_n\right)
    \leq p_{\!_A}  (1-X_n)+\rho_{n+1}^2(1-p_{\!_B}).
    \]
\end{Pro}
\noindent {\sc Proof:} We have
\[
\D M_{n+1}=V_{n+1}-\E(V_{n+1}\;|\;\F_n)+W_{n+1}-\E(W_{n+1}\;|\;\F_n),
\]
with
\[
V_{n+1}=\ind_{\{U_{n+1}\leq X_n\}\cap A_{n+1}}(1-X_n)
          -\ind_{\{U_{n+1}> X_n\}\cap B_{n+1}}X_n
	\]
	and
\[
W_{n+1}=-\rho_{n+1}\left(
	X_n\ind_{\{U_{n+1}\leq X_n\}\cap A^c_{n+1}}
	-(1-X_n)\ind_{\{U_{n+1}>X_n\}\cap B^c_{n+1}}\right).
\]
Note that $V_{n+1}W_{n+1}=0$, so that
\begin{eqnarray*}
 \E\left(\D M_{n+1}^2\;|\;\F_n\right)&=&\E(V_{n+1}^2\;|\;\F_n)+\E(W_{n+1}^2\;|\;\F_n)
 -\left(\E(V_{n+1}+W_{n+1}\;|\;\F_n)\right)^2\\
 &\leq &\E(V_{n+1}^2\;|\;\F_n)+\E(W_{n+1}^2\;|\;\F_n).
\end{eqnarray*}
Now, using $p_{\!_B}\leq p_{\!_A}$ and $X_n\leq 1$,
\begin{eqnarray*}
 \E\left(V_{n+1}^2\;|\;\F_n\right)&=&p_{\!_A}X_n(1-X_n)^2+p_{\!_B}(1-X_n)X_n^2\\
 &\leq &p_{\!_A}(1-X_n)\\
\mbox{ and }\hskip 2 cm     \E(W_{n+1}^2\;|\;\F_n)
&=&\rho_{n+1}^2\left(X_n^3(1-p_{\!_A})+(1-X_n)^3(1-p_{\!_B})\right)\hskip
2,25 cm \\
 &\leq&\rho_{n+1}^2(1-p_{\!_B}).
\end{eqnarray*}
This proves the Proposition.$\cqfd$
\section{The rate of convergence: pointwise 
convergence}
\label{SectionPtw}
\subsection{Convergence in probability}
\begin{Thm}\label{Conv-in-Prob}
    Assume 
    \begin{equation}\label{rho-conditions}
    \lim_{n\to\infty}\rho_n=0, \quad\lim_{n\to\infty}\frac{\gamma_n}{\rho_n}=0,\quad 
    \sum_n\rho_n\gamma_n=\infty, \quad\rho_n-\rho_{n-1}=o(\rho_n\gamma_n).
    \end{equation}
    Then, the sequence $((1-X_n)/\rho_n)_{n\geq 1}$ converges to 
    $(1-p_{\!_A})/\pi$ in probability.
\end{Thm}
Note that the assumptions of Theorem~\ref{Conv-in-Prob}
are satisfied if $\gamma_n=C/n^a$ and $\rho_n=C'/n^r$,
with $C,C'>0$, $0<r<a$ and $a+r<1$. In fact, we will see
that for this choice of parameters, convergence holds with probability
one (see Theorem~\ref{almost-sure}).

Before proving  Theorem~\ref{Conv-in-Prob}, we introduce the 
notation
 \[
 Y_n=\frac{1-X_n}{\rho_n}.
  \]
  We have, from (\ref{mean-algo})
\begin{eqnarray*}
    1-X_{n+1}&=&1-X_n-\gamma_{n+1}\pi X_n(1-X_n)-
       \rho_{n+1}\gamma_{n+1}\kappa(X_n)-\gamma_{n+1}\D M_{n+1}\\
       \frac{1-X_{n+1}}{\rho_{n+1}}&=&\frac{1-X_n}{\rho_{n+1}}
       -\frac{\gamma_{n+1}}{\rho_{n+1}}\pi X_n(1-X_n)
      - \gamma_{n+1}\kappa(X_n)-\frac{\gamma_{n+1}}{\rho_{n+1}}\D M_{n+1}.
\end{eqnarray*}
Hence
\begin{eqnarray*}
    Y_{n+1}&=&Y_n+(1-X_n)\left(\frac{1}{\rho_{n+1}}-
       \frac{1}{\rho_n}-\frac{\gamma_{n+1}}{\rho_{n+1}}\pi X_n\right)-
       \gamma_{n+1}\kappa(X_n)-\frac{\gamma_{n+1}}{\rho_{n+1}}\D M_{n+1}\\
       Y_{n+1}&=&Y_n\left(1+\gamma_{n+1}\e_n-\pi_n\gamma_{n+1} X_n\right)-
     \gamma_{n+1}\kappa(X_n)-\frac{\gamma_{n+1}}{\rho_{n+1}}\D M_{n+1},
\end{eqnarray*}
where 
\[
\e_n=\frac{\rho_n}{\gamma_{n+1}}\left(\frac{1}{\rho_{n+1}}-
       \frac{1}{\rho_n}\right)\mbox{ and } 
       \pi_n=\frac{\rho_n}{\rho_{n+1}}\pi.
\]
It follows from the assumption $\rho_n-\rho_{n-1}=o(\rho_n\gamma_n)$ 
that $\displaystyle \lim_{n\to\infty}\e_n=0$ and
$\displaystyle \lim_{n\to\infty}\pi_n=\pi$.

\begin{Lem}\label{Lem-stopY}
    Consider  two positive numbers $\pi^-$ and $\pi^+$  
with $0<\pi^-<\pi<\pi^+<1$. Given $l\in\N$, let
\[
\nu^l=\inf\{n\geq l\;|\; \pi_nX_n-\e_n>\pi^+\mbox{ or } 
                 \pi_nX_n-\e_n<\pi^-\}.
\]
We have 
\begin{itemize}
    \item $\displaystyle \lim_{l\to\infty}\P(\nu^l=\infty)=1$,
    \item for $n\geq l$, if $\theta_n^+=\prod_{k=l+1}^n(1-\pi^+\gamma_k)$ 
and $\theta_n^-=\prod_{k=l+1}^n(1-\pi^-\gamma_k)$,
    \begin{equation}\label{Ytheta-}
\frac{Y_{n\wedge\nu^l}}{\theta_{n\wedge\nu^l}^-}
\leq Y_l-\sum_{k=l+1}^{n\wedge\nu^l}
  \frac{\gamma_{k}}{\theta_{k}^-}\kappa(X_{k-1})
  -\sum_{k=l+1}^{n\wedge\nu^l}
  \frac{\gamma_k}{\rho_k\theta_k^-}\D M_{k}
\end{equation}
and
\begin{equation}\label{Ytheta+}
\frac{Y_{n\wedge\nu^l}}{\theta_{n\wedge\nu^l}^+}
\geq Y_l-\sum_{k=l+1}^{n\wedge\nu^l}
  \frac{\gamma_{k}}{\theta_{k}^+}\kappa(X_{k-1})
  -\sum_{k=l+1}^{n\wedge\nu^l}
  \frac{\gamma_k}{\rho_k\theta_k^+}\D M_{k}.
\end{equation}
\end{itemize}
Moreover, with the notation $||k||_\infty=\sup_{0<x<1}|\kappa(x)|$,
\[
\sup_{n\geq l}\E\left(Y_n\ind_{\{\nu^{l}=\infty\}}\right)
      \leq \E \,Y_l+\frac{||k||_\infty}{\pi^-}.
\]
\end{Lem}
\begin{Rem}\label{RemTight}
    Note that, as the proof will show, Lemma~\ref{Lem-stopY}
    remains valid if the condition 
    $\displaystyle \lim_{n\to\infty}\gamma_n/\rho_n=0$
    in (\ref{rho-conditions}) is replaced by the boundedness of 
    the sequence $(\gamma_n/\rho_n)_{n\geq 1}$. In particular,
    the last statement, which implies the tightness of the sequence
    $(Y_n)_{n\geq 1}$, will be used in Section~\ref{Weak}.
\end{Rem}
\noindent{\sc Proof:}
Since $\displaystyle \lim_{n\to\infty}(\pi_nX_n-\e_n)=\pi$ a.s., we clearly have 
$\displaystyle \lim_{l\to\infty}\P(\nu^l=\infty)=1$. 

On the other hand,
for $l\leq n<\nu^l$, we have
\[
Y_{n+1}\leq Y_n(1-\gamma_{n+1}\pi^-)-\gamma_{n+1}\kappa(X_n)-
   \frac{\gamma_{n+1}}{\rho_{n+1}}\D M_{n+1}
\]
and
\[
Y_{n+1}\geq Y_n(1-\gamma_{n+1}\pi^+)-\gamma_{n+1}\kappa(X_n)-
   \frac{\gamma_{n+1}}{\rho_{n+1}}\D M_{n+1},
\]
so that, with the notation 
$\theta_n^+=\prod_{k=l+1}^n(1-\pi^+\gamma_k)$ 
and $\theta_n^-=\prod_{k=l+1}^n(1-\pi^-\gamma_k)$,
\[
\frac{Y_{n+1}}{\theta_{n+1}^-}\leq 
\frac{Y_n}{\theta_n^-}-\frac{\gamma_{n+1}}{\theta_{n+1}^-}\kappa(X_n)-
   \frac{\gamma_{n+1}}{\rho_{n+1}\theta_{n+1}^-}\D M_{n+1}
\]
and
\[
\frac{Y_{n+1}}{\theta_{n+1}^+}\geq 
\frac{Y_n}{\theta_n^+}-\frac{\gamma_{n+1}}{\theta_{n+1}^+}\kappa(X_n)-
   \frac{\gamma_{n+1}}{\rho_{n+1}\theta_{n+1}^+}\D M_{n+1}.
\]
By summing up these inequalities, we get (\ref{Ytheta-})
and 
(\ref{Ytheta+}).

By taking expectations in (\ref{Ytheta-}), we get
\begin{eqnarray*}
    \E \frac{Y_{n\wedge\nu^l}}{\theta_{n\wedge\nu^l}^-}
    &\leq &\E \,Y_l+||k||_\infty\E\sum_{k=l+1}^{n\wedge\nu^l}
  \frac{\gamma_{k}}{\theta_{k}^-}\\
  &=&\E \,Y_l + \frac{||k||_\infty}{\pi^-}\E\sum_{k=l+1}^{n\wedge\nu^l}
  \left(\frac{1}{\theta_{k}^-}-\frac{1}{\theta_{k-1}^-}\right)\\
   &\leq& \E \,Y_l + \frac{||k||_\infty}{\pi^-}\frac{1}{\theta_n^-}.
\end{eqnarray*}
We then have
\begin{eqnarray*}
    \E(Y_n\ind_{\{\nu^l=\infty\}})=
    \theta_n^-\E\left(
             \frac{Y_{n\wedge\nu^l}}{\theta_{n\wedge\nu^l}^-}
	      \ind_{\{\nu^l=\infty\}}\right)
	      &\leq 
	      &\theta_n^-\E\frac{Y_{n\wedge\nu^l}}{\theta_{n\wedge\nu^l}^-}\\
	      &\leq &\theta_n^-\left(
	      \E \,Y_l + \frac{||k||_\infty}{\pi^-}\frac{1}{\theta_n^-}
	              \right)\\
		    &\leq & \E \,Y_l+\frac{||k||_\infty}{\pi^-}.
    \end{eqnarray*}
\hfill$\cqfd$
\begin{Lem}\label{LemGaThRho}
    Let $(\theta_n)_{n\in\N}$ be a sequence of positive numbers such 
    that $\theta_n=\prod_{k=1}^n(1-p\gamma_k)$ for some $p\in(0,1)$.
    The sequence $\left(\theta_n\sum_{k=1}^n 
    \frac{\gamma_k}{\theta_k\rho_k}\D M_k\right)_{n\in\N}$ converges 
    to 0 in probability.
\end{Lem}
\noindent{\sc Proof:} It suffices to show convergence to 0 in probability 
for the associated conditional variances $T_n$, defined by
\[
T_n=\theta_n^2\sum_{k=1}^n
   \frac{\gamma_k^2}{\theta_k^2\rho_k^2}
   \E\left(\D M_k^2\;|\;\F_{k-1}\right).
\]
We know from  Proposition~\ref{ProVarCond} that
\begin{eqnarray*}
 \E\left(\D M_k^2\;|\;\F_{k-1}\right)&\leq &
 p_{\!_A}  (1-X_{k-1})+\rho_{k}^2(1-p_{\!_B})\\
 &=&p_{\!_A}\rho_{k-1}Y_{k-1}+\rho_{k}^2(1-p_{\!_B}).
\end{eqnarray*}
Therefore, $T_n\leq p_{\!_A} T^{(1)}_n+(1-p_{\!_B})T^{(2)}_n$, where
\[
T^{(1)}_n=\theta_n^2\sum_{k=1}^n
   \frac{\gamma_k^2}{\theta_k^2\rho_k^2}\rho_{k-1}Y_{k-1}
\]
and 
\[
T^{(2)}_n=\theta_n^2\sum_{k=1}^n
   \frac{\gamma_k^2}{\theta_k^2}.
\]
We first prove that $\displaystyle \lim_{n\to \infty}T^{(2)}_n=0$.
Note that, since $p\gamma_k\leq 1$, 
\begin{equation}\label{theta2}
\frac{1}{\theta_k^2}-\frac{1}{\theta_{k-1}^2}=
\frac{2p\gamma_k-p^2\gamma_k^2}{\theta_k^2}\geq 
p\frac{\gamma_k}{\theta_k^2}.
\end{equation}
Therefore,
\[
T^{(2)}_n\leq \frac{\theta_n^2}{p}\sum_{k=1}^n\gamma_k
      \left(\frac{1}{\theta_k^2}-\frac{1}{\theta_{k-1}^2}
      \right),
\]
and $\displaystyle \lim_{n\to \infty}T^{(2)}_n=0$ follows from Cesaro's lemma.

We now deal with $T^{(1)}_n$. First note that the assumption
$\rho_n-\rho_{n-1}=o(\rho_n\gamma_n)$ implies 
$\displaystyle \lim_{n\to\infty}\rho_n/\rho_{n-1}=1$, so that, the sequence 
$(\gamma_n)_{n\geq 1}$ being non-increasing with limit 0,
we only need to prove that 
$\displaystyle \lim_{n\to\infty}\bar{T}^{(1)}_n=0$ in probability,
where
\[
\bar{T}^{(1)}_n=\theta_n^2\sum_{k=1}^n
   \frac{\gamma_k^2}{\theta_k^2\rho_k}Y_k.
   \]
   Now, with the notation of Lemma~\ref{Lem-stopY}, we have, for
   $n\geq l>1$ and $\e>0$,
   \begin{eqnarray*}
       \P\left(\bar{T}^{(1)}_n\geq \e\right)
       &\leq &
       \P(\nu^l<\infty)+\P\left(
       \theta_n^2\sum_{k=1}^n
   \frac{\gamma_k^2}{\theta_k^2\rho_k}Y_k\ind_{\{\nu^l=\infty\}}\geq 
   \e\right)\\
   &\leq &
      \P(\nu^l<\infty)
      +\frac{1}{\e}\theta_n^2\sum_{k=1}^n
  \frac{\gamma_k^2}{\theta_k^2\rho_k}
  \E\left(Y_k\ind_{\{\nu^l=\infty\}}\right).
   \end{eqnarray*}
Using Lemma~\ref{Lem-stopY}, $\displaystyle \lim_{n\to\infty}{\gamma_n/\rho_n}=0$ and (\ref{theta2}),
we have
\[
\lim_{n\to\infty}\theta_n^2\sum_{k=1}^n
  \frac{\gamma_k^2}{\theta_k^2\rho_k}
  \E\left(Y_k\ind_{\{\nu^l=\infty\}}\right)=0.
\]
We also know that $\displaystyle \lim_{l\to\infty} \P(\nu^l<\infty)=0$. Hence,
\[
\hskip 5 cm \lim_{n\to\infty}\P\left(\bar{T}^{(1)}_n\geq
\e\right)=0.\hskip 5 cm
\cqfd
\]

\noindent{\sc Proof of Theorem \ref{Conv-in-Prob}:} First note that
if $\theta_n=\prod_{k=1}^n(1-p\gamma_k)$, with $0<p<1$, we have
\[
    \sum_{k=1}^n \frac{\gamma_k}{\theta_k}\kappa(X_{k-1})=
    \frac{1}{p}\sum_{k=1}^n 
    \left(
    \frac{1}{\theta_k}-\frac{1}{\theta_{k-1}}
    \right)\kappa(X_{k-1}).
\]
Hence, using $\displaystyle \lim_{n\to\infty}X_n=1$ and $\kappa(1)=-(1-p_{\!_A})$,
\[
\lim_{n\to\infty}\theta_n\sum_{k=1}^n 
\frac{\gamma_k}{\theta_k}\kappa(X_{k-1})=-\frac{1-p_{\!_A}}{p}.
\]
Going back to (\ref{Ytheta-}) and (\ref{Ytheta+}) and using Lemma~\ref{LemGaThRho}
with $p=\pi^+$ and $\pi^-$, 
and the fact that $\displaystyle \lim_{l\to\infty}\P(\nu^l=\infty)=1$, we have, for 
all $\e>0$, $\displaystyle \lim_{n\to\infty}\P(Y_n\geq \frac{1-p_{\!_A}}{\pi^-}+\e)=
\lim_{n\to\infty}\P(Y_n\leq \frac{1-p_{\!_A}}{\pi^+}-\e)=0$, and since
$\pi^+$ and $\pi^-$ can be made arbitrarily  close to $\pi$,
the Theorem is proved.\hfill$\cqfd$
\subsection{Almost sure convergence}
\begin{Thm}\label{almost-sure}
    In addition to {}(\ref{rho-conditions}), we assume that
    for all $\beta\in [0,1]$,
    \begin{equation}\label{beta}
\gamma_n\rho_n^\beta-\gamma_{n-1}\rho_{n-1}^\beta
=o(\gamma_n^2\rho_n^\beta),
    \end{equation}
    and that, for some $\eta>0$, we have
    \begin{equation}\label{eta}
    \forall C>0,\quad 
    \sum_n\exp\left(-C\frac{\rho_n^{1+\eta}}{\gamma_n}
                \right)<\infty.
    \end{equation}
    Then, with probability 1,
    \[
    \lim_{n\to\infty}\frac{1-X_n}{\rho_n}=\frac{1-p_{\!_A}}{\pi}.
    \]
\end{Thm}
Note that the assumptions of Theorem~\ref{almost-sure}
are satisfied if $\gamma_n=Cn^{-a}$ and $\rho_n=C'n^{-r}$,
with $C,C'>0$, $0<r<a$ and $a+r<1$.

The proof of Theorem~\ref{almost-sure} is based on the following lemma,
which will be proved later.
\begin{Lem}\label{iter}
    Under the assumptions of Theorem~\ref{almost-sure}, let 
    $\alpha\in [0,1]$ and let $(\theta_n)_{n\in\geq 1}$ be a sequence of 
    positive numbers such that $\theta_n=\prod_{k=1}^n(1-p\gamma_k)$,
    for some $p\in(0,1)$. On the set 
    $\{\sup_n(\rho_n^\alpha Y_n)<\infty\}$, we have
    \[
    \lim_{n\to\infty}\theta_n\rho_n^{\frac{\alpha-\eta}{2}-1}
       \sum_{k=1}^n\frac{\gamma_k}{\theta_k}\D M_k=0\mbox{ a.s.,}
    \]
    where $\eta$ satisfies (\ref{eta}).
\end{Lem}
\noindent{\sc Proof of Theorem~\ref{almost-sure}:}
We start from the following form of (\ref{mean-algo}):
\[
1-X_{n+1}=(1-X_n)(1-\gamma_{n+1}\pi 
         X_n)-\rho_{n+1}\gamma_{n+1}\kappa(X_n)-\gamma_{n+1}\Delta M_{n+1}.
\]
We know that $\displaystyle \lim_{n\to\infty}X_n=1$ a.s.. Therefore, given $\pi^+$
and $\pi^-$, with $0<\pi^-<\pi<\pi^+<1$, there exists $l\in \N$ such 
that, for $n\geq l$,
\[
1-X_{n+1}\leq (1-X_n)(1-\gamma_{n+1}\pi^-)
      -\rho_{n+1}\gamma_{n+1}\kappa(X_n)-\gamma_{n+1}\Delta M_{n+1}
\]
and
\[
1-X_{n+1}\geq (1-X_n)(1-\gamma_{n+1}\pi^+)
      -\rho_{n+1}\gamma_{n+1}\kappa(X_n)-\gamma_{n+1}\Delta M_{n+1},
\]
so that, with the notation 
$\theta_n^+=\prod_{k=l+1}^n(1-\pi^+\gamma_k)$ 
and $\theta_n^-=\prod_{k=l+1}^n(1-\pi^-\gamma_k)$,
\[
\frac{1-X_{n+1}}{\theta_{n+1}^-}\leq \frac{1-X_n}{\theta_n^-}
   -\frac{\rho_{n+1}\gamma_{n+1}}{\theta_{n+1}^-}\kappa(X_n)-
   \frac{\gamma_{n+1}}{\theta_{n+1}^-}\Delta M_{n+1}
\]
and
\[
\frac{1-X_{n+1}}{\theta_{n+1}^+}\geq \frac{1-X_n}{\theta_n^+}
   -\frac{\rho_{n+1}\gamma_{n+1}}{\theta_{n+1}^+}\kappa(X_n)-
   \frac{\gamma_{n+1}}{\theta_{n+1}^+}\Delta M_{n+1}.
\]
By summing up these inequalities, we get, for $n\geq l+1$,
\[
\frac{1-X_{n}}{\theta_{n}^-}\leq (1-X_l)
    -\sum_{k=l+1}^n \frac{\rho_k\gamma_k}{\theta_k^-}\kappa(X_{k-1})
    -\sum_{k=l+1}^n \frac{\gamma_{k}}{\theta_{k}^-}\Delta M_{k}
\]
and
\[
\frac{1-X_{n}}{\theta_{n}^+}\geq (1-X_l)
    -\sum_{k=l+1}^n \frac{\rho_k\gamma_k}{\theta_k^+}\kappa(X_{k-1})
    -\sum_{k=l+1}^n \frac{\gamma_{k}}{\theta_{k}^+}\Delta M_{k}.
\]
Hence
\begin{equation}\label{Ytheta-bis}
Y_{n}
\leq \frac{\theta_{n}^-}{\rho_n}(1-X_l)-\frac{\theta_n^-}{\rho_n}\sum_{k=l+1}^{n}
  \frac{\rho_k\gamma_{k}}{\theta_{k}^-}\kappa(X_{k-1})
  -\frac{\theta_n^-}{\rho_n}\sum_{k=l+1}^{n}
  \frac{\gamma_k}{\theta_k^-}\D M_{k}
\end{equation}
and
\begin{equation}\label{Ytheta+bis}
Y_{n}
\geq \frac{\theta_{n}^+}{\rho_n}(1-X_l)-\frac{\theta_n^+}{\rho_n}\sum_{k=l+1}^{n}
  \frac{\rho_k\gamma_{k}}{\theta_{k}^+}\kappa(X_{k-1})
  -\frac{\theta_n^+}{\rho_n}\sum_{k=l+1}^{n}
  \frac{\gamma_k}{\theta_k^+}\D M_{k}.
\end{equation}
We have, with probability 1, 
$\displaystyle \lim_{n\to\infty}\kappa(X_n)=\kappa(1)=-(1-p_{\!_A})$, and,
since $\sum_{n=1}^\infty \rho_n\gamma_n=+\infty$,
\begin{equation}\label{equiv1}
\sum_{k=l+1}^{n}
  \frac{\rho_k\gamma_{k}}{\theta_{k}^-}\kappa(X_{k-1})\sim -(1-p_{\!_A})
  \sum_{k=l+1}^{n}
  \frac{\rho_k\gamma_{k}}{\theta_{k}^-}.
\end{equation}
On the other hand,
\begin{eqnarray}
  \sum_{k=l+1}^{n}
  \frac{\rho_k\gamma_{k}}{\theta_{k}^-}&=&
    \frac{1}{\pi^-}\sum_{k=l+1}^{n}
    \rho_k\left(
  \frac{1}{\theta_{k}^-}-\frac{1}{\theta_{k-1}^-}\right)\nonumber\\
  &=&\frac{1}{\pi^-}\left(\sum_{k=l+1}^{n}(\rho_{k-1}-\rho_k)
              \frac{1}{\theta_{k-1}^-}+\frac{\rho_n}{\theta_n^-}
	               -\frac{\rho_l}{\theta_l^-}\right)\nonumber\\
		     &\sim&\frac{1}{\pi^-}\frac{\rho_n}{\theta_n^-}\label{equiv2},
\end{eqnarray}
where we have used the condition
$\rho_k\!-\!\rho_{k-1}\!=\!o(\rho_k\gamma_k)$. We deduce
from~(\ref{equiv1}) and~(\ref{equiv2}) that
\[
\lim_{n\to\infty}\frac{\theta_n^-}{\rho_n}\sum_{k=1}^n 
\frac{\rho_k\gamma_k}{\theta_k^-}\kappa(X_{k-1})=-\frac{1-p_{\!_A}}{\pi^-}
\]
and, also, that $\displaystyle \lim_{n\to\infty}\frac{\theta_{n}^-}{\rho_n}=0$.
By a similar argument, we get $\displaystyle \lim_{n\to\infty}\frac{\theta_{n}^+}{\rho_n}=0$ 
and
\[
\lim_{n\to\infty}\frac{\theta_n^+}{\rho_n}\sum_{k=1}^n 
\frac{\rho_k\gamma_k}{\theta_k^+}\kappa(X_{k-1})=-\frac{1-p_{\!_A}}{\pi^+}
\]
It follows from Lemma~\ref{iter}, that given $\alpha\!\in[0,1]$,
we have, on the set $E_\alpha:=\{\sup_n\!(\rho_n^\alpha
Y_n)\!<\!\infty\}$,
\[
\lim_{n\to\infty}\rho_n^{\frac{\alpha-\eta}{2}-1}\theta_n^\pm\sum_{k=1}^n 
\frac{\gamma_k}{\theta_k^\pm}\D M_k=0.
\]
Together with (\ref{Ytheta-bis}) and (\ref{Ytheta+bis}) this implies
\begin{itemize}
    \item
    $\displaystyle \lim_{n\to \infty} Y_n=(1-p_{\!_A})/\pi$ a.s., if $\frac{\alpha-\eta}{2}\leq 0$,
    \item
    $\displaystyle \lim_{n\to \infty} Y_n\rho_n^{\frac{\alpha-\eta}{2}}=0$ a.s., if 
    $\frac{\alpha-\eta}{2}> 0$.
\end{itemize}
We obviously have  $\P(E_\alpha)=1$ for $\alpha =1$. We deduce from
the previous argument that if $\P(E_\alpha)=1$ and $\frac{\alpha-\eta}{2}> 0$,
then $\P(E_{\alpha'})=1$, with $\alpha'=\frac{\alpha-\eta}{2}-1$. Set $\alpha_0=1$ and
$\alpha_{k+1}=\frac{\alpha_k-\eta}{2}-1$. If $\alpha_0\leq \eta$, we
have $\displaystyle \lim_{n\to \infty} Y_n=(1-p_{\!_A})/\pi$ a.s. on $E_{\alpha_0}$. If 
$\alpha_0> \eta$, let $j$ be the largest 
integer such that $\alpha_j>\eta$ (note that $j$ exists because 
$\displaystyle \lim_{k\to\infty}\alpha_k<0$). We have $\P(E_{\alpha_{j+1}})=1$, and,
on $E_{\alpha_{j+1}}$, $\displaystyle \lim_{n\to \infty} Y_n=(1-p_{\!_A})/\pi$ a.s., because 
$\frac{\alpha_j-\eta}{2}\leq 0$.\hfill$\cqfd$

\bigskip

We now turn to the proof of Lemma~\ref{iter} which is based on the
following classical martingale inequality (see~\cite{MASS}, remark~1, p.14 for a proof in the case of i.i.d.
random variables: the extension to bounded martingale increments is straightforward).
\begin{Lem} (Bernstein's inequality for bounded martingale increments)
    \label{Bernstein}
    Let $(Z_i)_{1\leq i\leq n}$ be a finite sequence of square integrable 
    random variables, adapted to the filtration $(\F_i)_{1\leq i\leq 
    n}$, such that
    \begin{enumerate}
        \item
        $\E(Z_i\;|\;\F_{i-1})=0$, $i=1,\ldots,n$,
        \item
        $\E(Z_i^2\;|\;\F_{i-1})\leq \sigma_i^2$, $i=1,\ldots,n$,
        \item
        $|Z_i|\leq \D_n$, $i=1,\ldots,n$,
    \end{enumerate}
    where $\sigma_1^2$, \ldots, $\sigma_n^2$, $\D_n$ are deterministic 
    positive constants.
    
    Then, the following inequality holds:
    \[
    \forall \lambda >0,\quad \P\left(\left|\sum_{i=1}^nZ_i\right|
    \geq \lambda\right) \leq 
    2\exp\left(-\frac{\lambda^2}{2\left(b_n^2+\lambda\frac{\D_n}{3}\right)}
          \right),
	\]
	with $b_n^2=\sum_{i=1}^n\sigma_i^2$.
 \end{Lem}
 We will also need the following technical result.
 \begin{Lem}\label{tech}
     Let $(\theta_n)_{n\geq 1}$ be a sequence of 
    positive numbers such that $\theta_n=\prod_{k=1}^n(1-p\gamma_k)$,
    for some $p\in(0,1)$ and let $(\xi_n)_{n\geq 1}$ be a sequence of 
    non-negative numbers satisfying
    \[
    \gamma_n\xi_n-\gamma_{n-1}\xi_{n-1}=o(\gamma_n^2\xi_n).
    \]
    We have
    \[
    \sum_{k=1}^n\frac{\gamma_k^2\xi_k}{\theta_k^2}\sim 
      \frac{\gamma_n\xi_n}{2p \theta_n^2}.
    \]
 \end{Lem}
 \noindent{\sc Proof:} First observe that the condition 
 $\gamma_n\xi_n-\gamma_{n-1}\xi_{n-1}=o(\gamma_n^2\xi_n)$
 implies $\gamma_n\xi_n\sim \gamma_{n-1}\xi_{n-1}$ and
 that, given $\e >0$, we have, for $n$ large enough,
  \begin{eqnarray*}
 \gamma_n\xi_n-\gamma_{n-1}\xi_{n-1}&\geq& -\e \gamma_n^2\xi_n\\
      &\geq & -\e\gamma_{n-1}\gamma_n\xi_n,
  \end{eqnarray*}
  where we have used the fact that the sequence $(\gamma_n)$
  is non-increasing.
Since $\gamma_n\xi_n\sim \gamma_{n-1}\xi_{n-1}$, we have, for $n$ 
large enough, say $n\geq n_0$,
 \[
     \gamma_n\xi_n\geq 
     \gamma_{n-1}\xi_{n-1}
     (1-2\e\gamma_{n-1}).
  \] 
 Therefore, for $n>n_0$,
 \[
 \gamma_n\xi_n\geq \gamma_{n_0}\xi_{n_0}\prod_{k=n_0+1}^n(1-2\e 
         \gamma_{k-1}).
 \]
 From this, we easily deduce that 
 $\displaystyle \lim_{n\to\infty}\gamma_n\xi_n/\theta_n=\infty$ and that
 $\sum_n \gamma_n^2\xi_n/\theta_n^2=\infty$.

 Now, from
 \[
 \frac{1}{\theta_{k}^2}- \frac{1}{\theta_{k-1}^2}=
 \frac{2\gamma_k p-\gamma_k^2p^2}{\theta_k^2},
 \]
 we deduce (recall that $\displaystyle \lim_{n\to \infty}\gamma_n=0$)
 \[
 \frac{\gamma_k^2\xi_k}{\theta_k^2}\sim 
   \frac{\gamma_k\xi_k}{2p}\left(\frac{1}{\theta_{k}^2}- \frac{1}{\theta_{k-1}^2}
                      \right),
\]
and, since $\sum_n \gamma_n^2\xi_n/\theta_n^2=\infty$,
\begin{eqnarray*}
\sum_{k=1}^n\frac{\gamma_k^2\xi_k}{\theta_k^2}&\sim &
      \frac{1}{2p}\sum_{k=1}^n \gamma_k\xi_k
          \left(\frac{1}{\theta_{k}^2}- \frac{1}{\theta_{k-1}^2}
                      \right)\\
		  &=&
		  \frac{1}{2p}\left(
		  \frac{\gamma_n\xi_n}{\theta_n^2}+
		 \sum_{k=1}^n (\gamma_{k-1}\xi_{k-1}
		        -\gamma_{k}\xi_{k})\frac{1}{\theta_{k-1}^2}
		  \right)\\
		  &=&
		  \frac{\gamma_n\xi_n}{2p\theta_n^2}+
		  o\left(\sum_{k=1}^n 
		        \frac{\gamma_{k}^2\xi_{k}}{\theta_{k}^2}
		  \right),
\end{eqnarray*}
where, for the first equality, we have assumed $\xi_0=0$, and, for the 
last one, we have used again $\gamma_n\xi_n-\gamma_{n-1}\xi_{n-1}=o(\gamma_n^2\xi_n)$.
 \hfill$\cqfd$
 
 \bigskip
 
\noindent{\sc Proof of Lemma~\ref{iter}:} Given $\mu >0$, let
\[
\nu_\mu=\inf\{ n\geq 0\;|\;\rho_n^\alpha Y_n>\mu\}.
\]
Note that $\{\sup_n \rho_n^\alpha 
Y_n<\infty\}=\bigcup_{\mu>0}\{\nu_\mu=\infty\}$.

On the set $\{\nu_\mu=\infty\}$, we have
\[
\sum_{k=1}^n \frac{\gamma_k}{\theta_k}\D M_k=
    \sum_{k=1}^n \frac{\gamma_k}{\theta_k}\ind_{\{k\leq\nu_\mu\}}\D M_k.
\]
We now apply  Lemma~\ref{Bernstein} with 
$Z_i=\frac{\gamma_i}{\theta_i}\ind_{\{i\leq\nu_\mu\}}\D M_i$.
We have, using Proposition~\ref{ProVarCond},
\begin{eqnarray*}
\E(Z_i^2\;|\;\F_{i-1})&=&\frac{\gamma_i^2}{\theta_i^2}\ind_{\{i\leq\nu_\mu\}}
   \E(\D M_i^2\;|\;\F_{i-1})  \\
   &\leq &\frac{\gamma_i^2}{\theta_i^2}\ind_{\{i\leq\nu_\mu\}}
   \left(p_{\!_A} \rho_{i-1}Y_{i-1}+\rho_i^2(1-p_{\!_B})\right)\\
   &\leq &
   \frac{\gamma_i^2}{\theta_i^2}
   \left(p_{\!_A} \rho_{i-1}^{1-\alpha}\mu+\rho_i^2(1-p_{\!_B})\right),
\end{eqnarray*}
where we have used the fact that, on $\{i\leq\nu_\mu\}$,
$\rho_{i-1}^\alpha Y_{i-1}\leq \mu$. Since $\displaystyle \lim_{n\to 
\infty}\rho_n=0$ and $\displaystyle \lim_{n\to\infty}\rho_n/\rho_{n-1}=1$
(which follows from $\rho_n-\rho_{n-1}=o(\gamma_n\rho_n)$), we
have
\[
\E(Z_i^2\;|\;\F_{i-1})\leq \sigma_i^2, 
\]
with $\sigma_i^2=C_\mu 
\frac{\gamma_i^2\rho_i^{1-\alpha}}{\theta_i^2}$,
for some $C_\mu>0$, depending only on $\mu$. Using Lemma~\ref{tech},
we have
\[
\sum_{i=1}^n\sigma_i^2\sim C_\mu \frac{\gamma_n\rho_n^{1-\alpha}}{2p\theta_n^2}.
\]
On the other hand, we have, because the jumps $\D M_i$ are bounded,
\[
|Z_i|\leq C\frac{\gamma_i}{\theta_i},
\]
for some $C>0$. Note that 
$\frac{\gamma_k/\theta_k}{\gamma_{k-1}/\theta_{k-1}}=
\frac{\gamma_k}{\gamma_{k-1}(1-p\gamma_k)}$, and, since 
$\gamma_k-\gamma_{k-1}=o(\gamma_k^2)$ (take $\beta =0$ in 
(\ref{beta})), we have, for $k$ large enough,
$\gamma_k-\gamma_{k-1}\geq -p\gamma_k\gamma_{k-1}$, so that
$\gamma_k/\gamma_{k-1}\geq 1-p\gamma_k$, and the sequence
$(\gamma_n/\theta_n)$ is non-increasing for $n$ large enough.
Therefore, we have
\[
\sup_{1\leq i\leq n}|Z_i|\leq \D_n,
\]
with $\D_n\!= \!C\gamma_n/\theta_n$ for some $C\!>\!0$.
Now,  applying Lemma~\ref{Bernstein} with 
$
\lambda\!=\!\lambda_0\rho_n^{1-\frac{\alpha-\eta}{2}}\!/\theta_n,
$
we get
\begin{eqnarray*}
    \P\left(
    \theta_n\left|
       \sum_{k=1}^n \frac{\gamma_k}{\theta_k}
           \ind_{\{k\leq \nu_\mu\}}\D M_k
	    \right|
	     \geq \lambda_0\rho_n^{1-\frac{\alpha-\eta}{2}}
       \right)
      & \leq&
       2\exp \left(
       -\frac{\lambda_0^2\rho_n^{2-\alpha+\eta}}{2\theta_n^2b_n^2+
           2\lambda_0 \theta_n\rho_n^{1-\frac{\alpha-\eta}{2}}\frac{\D_n}{3}}
	   \right)\\
	 &\leq&
	 2\exp\left(  
	 -\frac{C_1\rho_n^{2-\alpha+\eta}}{C_2\gamma_n\rho_n^{1-\alpha}+
           C_3 
           \gamma_n\rho_n^{1-\frac{\alpha-\eta}{2}}}
        \right)\\
        &\leq &
       2 \exp\left(-C_4\frac{\rho_n^{1+\eta}}{\gamma_n}\right),
\end{eqnarray*}
where the positive constants $C_1$, $C_2$, $C_3$ and $C_4$ depend
on $\lambda_0$ and $\mu$, but not on $n$. Using~(\ref{eta}) and the 
Borel-Cantelli lemma, we conclude that, on $\{\nu_\mu=\infty\}$,
we have, for $n$ large enough,
\[
\theta_n\left|
       \sum_{k=1}^n \frac{\gamma_k}{\theta_k}
           \D M_k
	    \right|
	     < \lambda_0\rho_n^{1-\frac{\alpha-\eta}{2}}, \mbox { a.s.,}
\]
and, since $\lambda_0$ is arbitrary,
this completes the proof of the Lemma.\hfill$\cqfd$
\section{Weak  convergence of the normalized algorithm}
\label{Weak}
Throughout this section, we assume
(in addition to the initial conditions on the sequence 
$(\gamma_n)_{n\in\N}$)
\begin{equation}
    \label{rho/gammacond}
\gamma_n^2-\gamma_{n-1}^2=o(\gamma_n^2)\quad\mbox{and}\quad
    \frac{\gamma_n}{\rho_n}=g+o(\gamma_n),
\end{equation}
    where $g$ is a positive constant. Note that
    a possible choice is $\gamma_n=ag/\sqrt{n}$ and 
    $\rho_n=a/\sqrt{n}$, with $a>0$.
    
    Under these conditions, we have 
    $\rho_n-\rho_{n-1}=o(\gamma_n^2)$, and  we 
    can write, as in the beginning of Section~\ref{SectionPtw},
    \begin{equation}\label{Y-dynamics}
    Y_{n+1}=Y_n\left(1+\gamma_{n+1}\e_n-\pi_n\gamma_{n+1} X_n\right)-
     \gamma_{n+1}\kappa(X_n)-\frac{\gamma_{n+1}}{\rho_{n+1}}\D M_{n+1},
     \end{equation}
where 
 $\displaystyle \lim_{n\to\infty}\e_n=0$ and $\displaystyle \lim_{n\to\infty}\pi_n=\pi$.
 As observed in Remark~\ref{RemTight}, we know that, under the 
 assumptions (\ref{rho/gammacond}), the sequence $(Y_n)_{n\geq 1}$
 is tight.
 We will prove that it is convergent in distribution.
 \begin{Thm}\label{TCL}
 Under conditions (\ref{rho/gammacond}), the sequence $(Y_n=(1-X_n)/\rho_n)_{n\in\N}$
 converges weakly to the unique stationary distribution  of the Markov 
 process on $[0,+\infty)$ with generator $L$ defined by
 \begin{equation}\label{generateur}
 Lf(y)=p_{\!_B}y\frac{f(y+g)-f(y)}{g}+(1-p_{\!_A}-p_{\!_A}y)f'(y), \quad y\geq 
        0,
 \end{equation}
 for $f$ continuously differentiable and compactly supported in 
 $[0,+\infty)$.
 \end{Thm}
 The method for proving Theorem~\ref{TCL} is 
 based on the classical  functional approach to central limit
 theorems for stochastic algorithms (see Bouton~\cite{BOU}, Kushner~\cite{KUYI},
Duflo~\cite{DUF}).
 The long time behavior of the sequence $(Y_n)$ will be elucidated
 through the study of a sequence of continuous-time processes
 $Y^{(n)}=(Y^{(n)}_t)_{t\geq 0}$, which will be proved to converge
 weakly to the Markov process with generator $L$.
 We will show that $\nu$ has a unique stationary distribution,
 and that this is the weak limit of the sequence $(Y_n)_{n\in\N}$.
 
 The sequence $Y^{(n)}$ is defined as follows. Given $n\in\N$,
 and $t\geq 0$, set
  \begin{equation}\label{Y(n)}
      Y^{(n)}_t=Y_{N(n,t)}, 
   \end{equation}
   where
 \[
 N(n,t)=\min\left\{m\geq n\;|\;\sum_{k=n}^m\gamma_{k+1}>t\right\},
 \]
 so that $N(n,0)=n$, for $t\in[0,\gamma_{n+1})$, and, for
 $m\geq n+1$, $N(n,t)=m$ if and only if $\sum_{k=n+1}^{m}\gamma_k\leq t
   <\sum_{k=n+1}^{m+1}\gamma_k$. 
   \begin{Thm}\label{convY(n)}
       Under the assumptions of Theorem~\ref{TCL}, the sequence  
       of continuous time processes $(Y^{(n)})_{n\in\N}$ converges
       weakly (in the sense of Skorokhod) to a Markov process with generator 
       $L$.
    \end{Thm}
    The proof of Theorem~\ref{convY(n)} is done in two steps: in 
    section~\ref{Tight}, we prove tightness,  in section~\ref{Identif}, we characterize the
    limit by a martingale problem.
 \subsection{Tightness}
 \label{Tight}
 It follows
 from (\ref{Y-dynamics}) that the process $Y^{(n)}$ admits the 
 following decomposition:
 \begin{equation}\label{Y-decomp}
     Y^{(n)}_t=Y_n+B^{(n)}_t+M^{(n)}_t,
 \end{equation}
 with
 \[
 B^{(n)}_t=-\sum_{k=n+1}^{N(n,t)}\gamma_k\left[
      \kappa(X_{k-1})+(\pi_{k-1}X_{k-1}-\e_{k-1})Y_{k-1}\right]
 \]
 and
 \[
 M^{(n)}_t=-\sum_{k=n+1}^{N(n,t)}\frac{\g_k}{\rho_k}\D M_k.
 \]
 The process $(M^{(n)}_t)_{t\geq 0}$ is a square integrable martingale with respect to 
 the filtration $(\F^{(n)}_t)_{t\geq 0}$, with 
 $\F^{(n)}_t=\F_{N(n,t)}$, and we have
 \[
 <\!\!M^{(n)}\!\!>_t=\sum_{k=n+1}^{N(n,t)}\left(\frac{\g_k}{\rho_k}\right)^2
      \E(\D M_k^2 \;|\; \F_{k-1}).
 \]
 We already know (see Remark~\ref{RemTight}) that the sequence 
 $(Y_n)_{n\in\N}$ is tight. 
 Recall that in order for the sequence $(M^{(n)})$ to be tight,
 it is sufficient that the sequence
 $(<\!\!M^{(n)}\!\!>)$ is $C$-tight (see \cite{JacodShi}, Theorem 
 4.13, p. 358, chapter VI). Therefore, the tightness
 of the sequence $(Y^{(n)})$  in the sense of Skorokhod will follow
 from the following result.
 \begin{Pro}\label{Pro-Tight}
     Under the assumptions (\ref{rho/gammacond}),  the sequences
      $(B^{(n)})$ and $(<\!\!M^{(n)}\!\!>)$ are $C$-tight.
 \end{Pro}
 For the proof of this proposition,we will need the following lemma.
 \begin{Lem}\label{Lem-Tightness}
     Define $\nu^l$ as in Lemma~\ref{Lem-stopY}, for $l\in\N$.
              There exists a positive constant $C$ such that, for all
          $l, n, N\in \N$ with  $l\leq n\leq N$, we have
	\[
	\forall \lambda \geq 1, \quad
	\P\left(\sup_{n\leq j\leq N}|Y_j-Y_n|\geq \lambda\right)
	\leq \P(\nu^l<\infty)+C
	\frac{(1+\E \,Y_l)\left(\sum_{k=n+1}^N\gamma_k\right)}{\lambda}.
	\]
 \end{Lem}
 \noindent{\sc Proof:} 
 The function $\kappa$ being bounded on $[0,1]$, it follows
  from (\ref{Y-dynamics}) that there exist positive, deterministic
  constants $a$ and $b$ such that, for all $n\in\N$,
 \begin{equation}\label{Y-control}
     -\gamma_{n+1}(a+bY_n)-\frac{\gamma_{n+1}}{\rho_{n+1}}\D M_{n+1}
     \leq
     Y_{n+1}-Y_n
     \leq
     \gamma_{n+1}(a+bY_n)-\frac{\gamma_{n+1}}{\rho_{n+1}}\D M_{n+1}.
 \end{equation}
 We also know from Proposition~\ref{ProVarCond} that
 \begin{equation}\label{Y-VarCond}
     \E\left( \D M_{n+1}^2\;|\;\F_n\right)
     \leq p_{\!_A}\rho_n Y_n+(1-p_{\!_B})\rho_{n+1}^2.
 \end{equation}

 From (\ref{Y-control}), we derive, for $j\geq n$,
 \[
 |Y_j-Y_n|\leq \sum_{k=n+1}^j\gamma_k (a+bY_{k-1})+
   \left|
   \sum_{k=n+1}^j\frac{\gamma_k}{\rho_k}\D M_k.
   \right|
 \]
 Let $\tilde{Y}_k=Y_k\ind_{\{k\leq \nu^l\}}$ and 
 $\D \tilde{M}_k=\ind_{\{k\leq \nu^l\}}\D M_k$
 On the set $\nu^l=\infty$, we have $Y_{k-1}=\tilde{Y}_{k-1}$
 and $\D M_k=\D \tilde{M}_k$. Hence
 \begin{eqnarray*}
 \P\left(\sup_{n\leq j\leq N}|Y_j-Y_n|\geq \lambda\right)
 &\leq &\P(\nu^l<\infty)+\P\left(\sum_{k=n+1}^N\gamma_k (a+b\tilde{Y}_{k-1})\geq
      \lambda/2\right)+\\
      &&
      \P\left(\sup_{n\leq j\leq N}\left|
      \sum_{k=n+1}^j\frac{\gamma_k}{\rho_k}\D \tilde{M}_k\right|
      \geq
      \lambda/2\right).
 \end{eqnarray*}
 We have, using Markov's inequality and Lemma~\ref{Lem-stopY},
 \begin{eqnarray*}
\P\left(\sum_{k=n+1}^N\gamma_k (a+b\tilde{Y}_{k-1})\geq
      \lambda/2\right)
      &\leq &\frac{2}{\lambda}\E\sum_{k=n+1}^N\gamma_k 
      (a+b\tilde{Y}_{k-1})\\
      &\leq &\frac{2}{\lambda}\left( a+b\sup_{k\geq l}\E 
      \left(Y_k\ind_{\{\nu^l=\infty\}}\right)\right)\sum_{k=n+1}^N\gamma_k\\
      &\leq&\frac{2}{\lambda}\left(b\E \,Y_l 
      +b\frac{||\kappa||_\infty}{\pi^-}+a\right)\sum_{k=n+1}^N\gamma_k.
 \end{eqnarray*}
 On the other hand, using Doob's inequality,
 \begin{eqnarray*}
 \P\left(\sup_{n\leq j\leq N}\left|
      \sum_{k=n+1}^j\frac{\gamma_k}{\rho_k}\D \tilde{M}_k\right|
      \geq
      \lambda/2\right)
      &\leq &\frac{16}{\lambda^2}\E\sum_{k=n+1}^N
        \frac{\gamma_k^2}{\rho_k^2} \E\left( \D 
        \tilde{M}_{k}^2\;|\;\F_{k-1}\right)\\
     &\leq &\frac{16}{\lambda^2}\E\sum_{k=n+1}^N
        \frac{\gamma_k^2}{\rho_k^2}
        \ind_{\{k\leq \nu\}}\left(  
         p_{\!_A}\rho_{k-1} Y_{k-1}+(1-p_{\!_B})\rho_k^2\right).
 \end{eqnarray*}
 Using $\displaystyle \lim_n(\gamma_n/\rho_n)=g$, 
 $\rho_{k-1}\sim \rho_k$, $\displaystyle \lim_n\rho_n=0$ and Lemma~\ref{Lem-stopY}, we get, for some $C>0$,
 \[
 \P\left(\sup_{n\leq j\leq N}\left|
      \sum_{k=n+1}^j\frac{\gamma_k}{\rho_k}\D \tilde{M}_k\right|
      \geq
      \lambda/2\right)
\leq C\frac{(1+\E \,Y_l)\left(\sum_{k=n+1}^N\gamma_k\right)}{\lambda^2},
\]
and, since we have assumed $\lambda\geq 1$, the proof of the lemma is completed. \hfill$\cqfd$

\medskip

\noindent{\sc Proof of Proposition~\ref{Pro-Tight}:} 
Given $s$ and $t$, with $0\leq s\leq t$, we have, using the 
boundedness of $\kappa$,
\[
|B^{(n)}_t-B^{(n)}_s|\leq \sum_{k=N(n,s)+1}^{N(n,t)} \gamma_k(a+bY_{k-1})
\]
for some $a,b >0$.

Similarly, using (\ref{Y-VarCond}), we have
\[
\left|<\!\!M^{(n)}\!\!>_t-<\!\!M^{(n)}\!\!>_s\right|
\leq \sum_{k=N(n,s)+1}^{N(n,t)} \gamma_k(a'+b'Y_{k-1})
\]
for some $a',b' >0$. These inequalities express the fact
that the processes
 $B^{(n)}$ and $<\!\!M^{(n)}\!\!>$ are {\em strongly dominated}
 (in the sense of \cite{JacodShi}, definition 3.34)
 by a linear combination of the processes 
 $X^{(n)}$ and $Z^{(n)}$, where 
 $X^{(n)}_t=\sum_{k=n+1}^{N(n,t)}\gamma_k$ and 
 $Z^{(n)}_t=\sum_{k=n+1}^{N(n,t)}\gamma_k Y_{k-1}$. Therefore, we 
 only need to prove that the sequences $(X^{(n)})$ and $(Z^{(n)})$
 are $C$-tight. This is obvious for the sequence $X^{(n)}$,
 which in fact converges to the deterministic process $t$.
 We now prove that $Z^{(n)}$ is $C$-tight. We have, for
 $0\leq s\leq t\leq T$
 \begin{eqnarray*}
 |Z^{(n)}_t-Z^{(n)}_s|&\leq& \left(\sup_{n\leq j\leq 
 N(n,T)}Y_j\right)\sum_{k=N(n,s)+1}^{N(n,t)}\gamma_k\\
 &\leq &(t-s+\gamma_{N(n,s)+1})\sup_{n\leq j\leq 
 N(n,T)}Y_j\\
 &\leq & (t-s+\gamma_{n+1})\sup_{n\leq j\leq 
 N(n,T)}Y_j,
 \end{eqnarray*}
where we have used $\sum_{k=n+1}^{N(n,t)}\gamma_k\leq t$ and 
$s\leq \sum_{k=n+1}^{N(n,s)+1}\gamma_k$ and the monotony of the 
sequence $(\gamma_n)_{n\geq 1}$.
 
 Therefore, for 
$\delta>0$, and $n$ large enough so that $\gamma_{n+1}\leq \delta$,
 \begin{eqnarray*}
    \P\left(\sup_{0\leq s\leq t\leq T, t-s\leq \delta} 
    |Z^{(n)}_t-Z^{(n)}_s|\geq    \eta\right)
        &\leq &\P\left(\sup_{n\leq j\leq N(n,T)}Y_j\geq 
	         \frac{\eta}{\delta +\gamma_{n+1}}\right)\\
       &\leq &
          \P\left(Y_n \geq 
	         \frac{\eta}{4\delta }\right)
		  \\
		  &&+
		  \P\left(\sup_{n\leq j\leq N(n,T)}|Y_j-Y_n|\geq 
				   \frac{\eta}{4\delta }\right).
 \end{eqnarray*}
We have, from Lemma~\ref{Lem-Tightness},
\begin{eqnarray*}
\P\left(\sup_{n\leq j\leq N(n,T)}|Y_j-Y_n|\geq 
				   \frac{\eta}{4\delta}\right)
				   &\leq&
				   \P(\nu^l<\infty)
				   +\frac{4C\delta}{\eta}(1+\E \,Y_l)\sum_{k=n+1}^{N(n,T)}\gamma_k\\
				   &\leq &
				   \P(\nu^l<\infty)+\frac{4CT\delta}{\eta}(1+\E \,Y_l).
\end{eqnarray*}	
We easily conclude from these estimates that, given $T>0$, $\e>0$ and 
$\eta>0$, we have for $n$ large enough and $\delta$ small enough,
\[
\P\left(\sup_{0\leq s\leq t\leq T, t-s\leq \delta} 
    |Z^{(n)}_t-Z^{(n)}_s|\geq    \eta\right)<\e,
\]
which proves the $C$-tightness of the sequence $(Z^{(n)})$.\hfill$\cqfd$
\subsection{Identification of the limit}
\label{Identif}
\begin{Lem}\label{Lem-Ident}
    Let $f$ be a $C^1$ function with compact support in $[0,+\infty)$.
    We have
    \[
    \E\left(f(Y_{n+1})-f(Y_n)\;|\;\F_n\right)=
       \gamma_{n+1}Lf(Y_n)+\gamma_{n+1}Z_n,\quad n\in\N,
    \]
    where the operator $L$ is defined by
    \begin{equation}\label{generator}
        Lf(y)=p_{\!_B}y\frac{f(y+g)-f(y)}{g}+(1-p_{\!_A}-p_{\!_A}y)f'(y), \quad y\geq 
        0,
    \end{equation}
    and the sequence $(Z_n)_{n\in\N}$ satisfies 
    $\displaystyle \lim_{n\to\infty}Z_n=0$ in probability.
\end{Lem}
\noindent {\sc Proof:} From (\ref{Y-dynamics}), we have
\begin{eqnarray}
    Y_{n+1}&=&Y_n+\gamma_{n+1}(-\kappa(1)-\pi 
    Y_n)-\frac{\gamma_{n+1}}{\rho_{n+1}}\D 
    M_{n+1}+\gamma_{n+1}\zeta_n\nonumber\\
    &=&Y_n+\gamma_{n+1}(1-p_{\!_A}-\pi 
    Y_n)-\frac{\gamma_{n+1}}{\rho_{n+1}}\D 
    M_{n+1}+\gamma_{n+1}\zeta_n\nonumber\\
    &=&Y_n+\gamma_{n+1}(1-p_{\!_A}-\pi 
    Y_n)-g\D 
    M_{n+1}+\gamma_{n+1}\zeta_n+
    \left(g-\frac{\gamma_{n+1}}{\rho_{n+1}}\right)\D 
    M_{n+1},\label{Y-dynamics1}
\end{eqnarray}
where $\zeta_n=\kappa(1)-\kappa(X_n)+Y_n(\pi-(\pi_n X_n-\e_n))$, so 
that $\zeta_n$ is $\F_n$-measurable and $\displaystyle \lim_{n\to\infty}\zeta_n=0$
in probability. Going back to (\ref{M-dynamics}), we rewrite the 
martingale increment $\D M_{n+1}$ as follows:
\begin{eqnarray*}
    \D M_{n+1}&=&-X_n\left(\ind_{\{U_{n+1}>X_n\}\cap 
        B_{n+1}}-p_{\!_B}(1-X_n)\right)
          +\rho_nY_n\left(\ind_{\{U_{n+1}\leq X_n\}\cap 
        A_{n+1}}-p_{\!_A} X_n\right)\\
        &&-\rho_{n+1}\left(X_n\ind_{\{U_{n+1}\leq X_n\}\cap 
        A_{n+1}^c}-(1-X_n)\ind_{\{U_{n+1}> X_n\}\cap B_{n+1}^c}
            +\kappa(X_n)\right).
\end{eqnarray*}
Hence,
\[
    Y_{n+1}=Y_n+\gamma_{n+1}(1-p_{\!_A}-\pi 
    Y_n+\zeta_n)+\xi_{n+1} +\D \hat{M}_{n+1},
\]
where
\[
\xi_{n+1}=gX_n\left(\ind_{\{U_{n+1}>X_n\}\cap 
        B_{n+1}}-p_{\!_B}(1-X_n)\right)
\]
and 
\begin{eqnarray*}
    \D \hat{M}_{n+1}&=&\left(g-\frac{\gamma_{n+1}}{\rho_{n+1}}\right)\D 
    M_{n+1}
          -g\rho_nY_n\left(\ind_{\{U_{n+1}\leq X_n\}\cap 
        A_{n+1}}-p_{\!_A} X_n\right)\\
        &&+g\rho_{n+1}\left(X_n\ind_{\{U_{n+1}\leq X_n\}\cap 
        A_{n+1}^c}-(1-X_n)\ind_{\{U_{n+1}> X_n\}\cap B_{n+1}^c}
            +\kappa(X_n)\right).
\end{eqnarray*}
Note that, due to our assumptions on $\gamma_n$ and $\rho_n$,
we have, for some deterministic positive constant $C$,
\begin{equation}\label{borne-M-hat}
    \left|\D \hat{M}_{n+1}\right|\leq C\gamma_{n+1}(1+Y_n), \quad 
    n\in \N.
\end{equation}
Now, let
\[
\tilde{Y}_n=Y_n+\gamma_{n+1}(1-p_{\!_A}-\pi 
    Y_n+\zeta_n)\mbox{ and } \bar{Y}_{n+1}=\tilde{Y}_n+\xi_{n+1},
\]
so that $Y_{n+1}=\bar{Y}_{n+1}+\D \hat{M}_{n+1}$.
We have
\[
f(Y_{n+1})-f(Y_n)=f(Y_{n+1})-f(\bar{Y}_{n+1})+f(\bar{Y}_{n+1})-f(Y_n).
\]
We will first show that 
\begin{equation}\label{bar1}
f(Y_{n+1})-f(\bar{Y}_{n+1})=f'(\tilde{Y}_n)\D 
\hat{M}_{n+1}+\gamma_{n+1}T_{n+1},\mbox{ where }
\P\mbox{-}\lim_{n\to\infty}\E(T_{n+1}\;|\;\F_n)=0, 
\end{equation} 
 with the notation $\P\mbox{-}\lim$ for a limit in probability.
Denote by $w$ the modulus of continuity of $f'$:
\[
w(\delta)=\sup_{|x-y|\leq \delta}|f'(y)-f'(x)|, \quad \delta>0.
\]
We have, for some (random) $\theta\in(0,1)$,
\begin{eqnarray*}
f(Y_{n+1})-f(\bar{Y}_{n+1})&=&f'(\bar{Y}_{n+1}+
        \theta \D \hat{M}_{n+1})\D \hat{M}_{n+1}\\
    &=&f'(\tilde{Y}_n)\D \hat{M}_{n+1}+ V_{n+1},
\end{eqnarray*}
where $V_{n+1}=\left(f'(\bar{Y}_{n+1}+
        \theta \D \hat{M}_{n+1})-f'(\tilde{Y}_n)\right)\D \hat{M}_{n+1}$.
We have
\begin{eqnarray*}
|V_{n+1}|&\leq &w\left(|\xi_{n+1}|+|\D \hat{M}_{n+1}|\right)|\D \hat{M}_{n+1}|\\
        &\leq &Cw\left(|\xi_{n+1}|+C\gamma_{n+1}(1+Y_n)\right)\gamma_{n+1}(1+Y_n),
\end{eqnarray*}
where we have used $\bar{Y}_{n+1}=\tilde{Y}_n+\xi_{n+1}$ and 
(\ref{borne-M-hat}). In order to get (\ref{bar1}), it suffices to 
prove that $\displaystyle \lim_{n\to\infty}
\E\left(w\left(|\xi_{n+1}|+C\gamma_{n+1}(1+Y_n)\right)\;|\;\F_n\right)=0$
in probability. On the set $\{U_{n+1}>X_n\}\cap B_{n+1}$, we have
$|\xi_{n+1}|=gX_n\left(1-p_{\!_B}(1-X_n)\right)\leq g$, and,  on the complement,
$|\xi_{n+1}|=gX_np_{\!_B}(1-X_n)\leq g(1-X_n)$. Hence
\begin{eqnarray*}
\E\left(w\left(|\xi_{n+1}|+C\gamma_{n+1}(1+Y_n)\right)\;|\;\F_n\right)
&\leq& p_{\!_B}(1-X_n)w\left(g+C\gamma_{n+1}(1+Y_n)\right)\\
&&
  +(1-p_{\!_B}(1-X_n))w\left(\hat{Y}_n\right),
\end{eqnarray*}
where $\hat{Y}_n=g(1-X_n)+C\gamma_{n+1}(1+Y_n)$. 
Observe that $\displaystyle \lim_{n\to\infty}\hat{Y}_n=0$ in  probability
(recall  that $\displaystyle \lim_{n\to\infty}X_n=1$ almost surely).
Therefore, we have (\ref{bar1}).

We deduce from 
$\E(\D \hat{M}_{n+1}\;|\;\F_n)=0$  that
\[
\E\left(f(Y_{n+1})-f(Y_n)\;|\;\F_n\right)=\gamma_{n+1}\E(T_{n+1}\;|\;\F_n)+
   \E\left(f(\bar{Y}_{n+1})-f(Y_n)\;|\;\F_n\right),
\]
so that the proof will be completed when we have shown
\begin{equation}
    \label{bar2}
   \P\mbox{-}\lim_{n\to\infty}\E\left(
  \frac{f(\bar{Y}_{n+1})-f(Y_n)-\gamma_{n+1}Lf(Y_n)}{\gamma_{n+1}}\;|\;\F_n\right)=0.
  \end{equation}
 We have
\begin{eqnarray*}
   \E\left(
  f(\bar{Y}_{n+1})\;|\;\F_n\right)&=&
     \E\left(
  f(\tilde{Y}_{n}+\xi_{n+1})\;|\;\F_n\right)\\
  &=&p_{\!_B}(1-X_n)f(\tilde{Y}_{n}+gX_n(1-p_{\!_B}(1-X_n)))\\
      *[.4em]&&+
      (1-p_{\!_B}(1-X_n))f(\tilde{Y}_{n}-gX_np_{\!_B}(1-X_n))\\
      &=&p_{\!_B}\rho_nY_nf(\tilde{Y}_{n}+gX_n(1-p_{\!_B}(1-X_n)))\\
      *[.4em]&&+
      (1-p_{\!_B}\rho_nY_n)f(\tilde{Y}_{n}-gX_np_{\!_B}(1-X_n)).
\end{eqnarray*}
Hence
\[
\E\left(
  f(\bar{Y}_{n+1})-f(Y_n)\;|\;\F_n\right)=F_n+G_n,
\]
with
\[
F_n=p_{\!_B}\rho_nY_n\left(f(\tilde{Y}_{n}+gX_n(1-p_{\!_B}(1-X_n)))-f(Y_n)\right)
\]
and
\[
G_n=(1-p_{\!_B}\rho_nY_n)\left(f(\tilde{Y}_{n}-gX_np_{\!_B}(1-X_n))
               -f(Y_n)\right).
\]
For the behavior of $F_n$ as $n$ goes to infinity, we use  
\[
\P\mbox{-}\lim_{n\to\infty}\left(\tilde{Y}_{n}+gX_n(1-p_{\!_B}(1-X_n))-Y_n-g\right)=0,
\]
and  $\displaystyle \lim_{n\to\infty}\rho_n/\gamma_{n+1}=1/g$, so 
that
\[
\P\mbox{-}\lim_{n\to\infty}\frac{F_n-p_{\!_B}Y_n\frac{f(Y_n+g)-f(Y_n)}{g}}{\gamma_{n+1}}=0.
\]
For the behavior of $G_n$, we write, using 
$\displaystyle \lim_{n\to\infty}\rho_n/\gamma_{n+1}=1/g$ again,
\begin{eqnarray*}
\tilde{Y}_{n}-gX_np_{\!_B}(1-X_n)&=&Y_n+\gamma_{n+1}\left(1-p_{\!_A}-\pi 
    Y_n+\zeta_n\right)-gp_{\!_B}X_n\rho_nY_n\\
    &=&Y_n+\gamma_{n+1}(1-p_{\!_A}-p_{\!_A}Y_n)+\gamma_{n+1}\eta_n,
\end{eqnarray*}
with $\P\mbox{-}\lim_{n\to\infty}\eta_n=0$, so that, using the fact
that $f$ is $C^1$ with compact support and the tightness of $(Y_n)$,
\[
\P\mbox{-}\lim_{n\to\infty}\frac{G_n-(1-p_{\!_A}-p_{\!_A}Y_n)f'(Y_n)}{\gamma_{n+1}}=0,
\]
which completes the proof of (\ref{bar2}).
\hfill$\cqfd$

\bigskip
\noindent{\sc Proof of Theorem~\ref{convY(n)}: }
As mentioned before, it follows from Proposition~\ref{Pro-Tight}
that the sequence of processes $(Y^{(n)})$ is tight in the Skorokhod 
sense. 

On the other hand, it follows from Lemma~\ref{Lem-Ident} that, if
     $f$ is a $C^1$ function with compact support in $[0,+\infty)$,
     we have
     \[
     f(Y_n)=f(Y_0)+\sum_{k=1}^n\gamma_k 
     Lf(Y_{k-1})+\sum_{k=1}^n\gamma_k Z_{k-1}+M_n,
     \]
     where $(M_n)$ is a martingale and $(Z_n)$ is an adapted sequence
     satisfying $\P$-$\displaystyle \lim_{n\to\infty}Z_n=0$. Therefore,
     \[
     f(Y^{(n)}_t)-f(Y^{(n)}_0)=M^{(n)}_t+\sum_{k=N(n,0)+1}^{N(n,t)}\gamma_k 
         (Lf(Y_{k-1})+Z_{k-1}),
     \]
     where $M^{(n)}_t=M_{N(n,t)}-M_{N(n,0)}$. It is easy to verify 
     that $M^{(n)}$ is a martingale with respect to $\F^{(n)}$.
     
     We also have
     \[
     \int_0^tLf(Y^{(n)}_s)ds=\sum_{k=n+1}^{N(n,t)}\gamma_k Lf(Y_{k-1})
     +\left(t-\sum_{k=n+1}^{N(n,t)}\gamma_k\right)f(Y^{(n)}_t).
     \]
     Therefore
     \[
     f(Y^{(n)}_t)-f(Y^{(n)}_0)-\int_0^tLf(Y^{(n)}_s)ds=M^{(n)}_t+R^{(n)}_t,
     \]
     where $\P$-$\displaystyle \lim_{n\to\infty}R^{(n)}_t=0$. It follows that any 
     weak limit of the sequence $(Y^{(n)})_{n\in \N}$ solves the 
     martingale problem associated with $L$. From this, together with 
     the study of the stationary distribution of $L$ (see 
     Section~\ref{Invariant}), we will deduce Theorem~\ref{TCL}
     and Theorem~\ref{convY(n)}.$\cqfd$
      \subsection{The stationary distribution}
      \label{Invariant}
      \begin{Thm}\label{stat-dist}
	The Markov process $(Y_t)_{t\geq 0}$, on $[0,+\infty)$, with generator $L$ has  a unique
	stationary probability distribution $\nu$. Moreover, 
	    $\nu$ has a density  on $[0,+\infty)$, which vanishes on 
	    $(0,r_{\!_A}]$ (where $r_{\!_A}=(1-p_{\!_A})/p_{\!_A}$), and is positive and 
	    continuous
	    on the open interval $(r_{\!_A},+\infty)$. The stationary 
	    distribution $\nu$ also satisfies the following property:
	   for every compact set $K$ in  $[0,+\infty)$, and every
	    bounded continuous function $f$, we have
	   \begin{equation}\label{unif}
	   \lim_{t\to\infty}\sup_{y\in K}\left|\E_y(f(Y_t))-\int f \,d\nu\right|=0.
	   \end{equation}
     \end{Thm}
     Before proving Theorem~\ref{stat-dist}, we will show how 
     Theorem~\ref{TCL} follows from~(\ref{unif}).
     
     \bigskip
     \noindent{\sc Proof of Theorem~\ref{TCL}:}
     Fix $t>0$. For $n$ large enough, we have 
     $\gamma_n\leq t <\sum_{k=1}^n\gamma_k$, so that there exists
     $\bar{n}\in \{1,\ldots,n-1\}$ such that
     \[
     \sum_{k=\bar{n}+1}^n\gamma_k\leq t<\sum_{k=\bar{n}}^n\gamma_k.
     \]
     Let $t_n=\sum_{k=\bar{n}+1}^n\gamma_k$. We have 
     \[
     0\leq t-t_n<\gamma_{\bar{n}}\quad\mbox{and}\quad 
           Y^{(\bar{n})}_{t_n}=Y_n.
     \]
     Since $t$ is fixed, the condition 
     $\sum_{k=\bar{n}+1}^n\gamma_k\leq t$ implies 
     $\displaystyle \lim_{n\to\infty}\bar{n}=\infty$ and $\displaystyle \lim_{n\to\infty}t_n=t$.
     
     Now, given $\e>0$, there is a compact set $K$ such that for 
     every weak limit $\mu$ of the sequence $(Y_n)_{n\in\N}$, 
     $\mu(K^c)<\e$. Using~(\ref{unif}),
     we choose $t$ such that 
     \[
     \sup_{y\in K}\left|\E_y(f(Y_t))-\int f d\nu\right|<\e.
     \]
     Now take a weakly convergent subsequence $(Y_{n_k})_{k\in\N}$.
     By another subsequence extraction, we can assume that the 
     sequence $(Y^{(\overline{n}_k)})$ converges weakly to a process $Y^{(\infty)}$
     which satisfies the martingale problem associated with $L$.
     We then have, due to the quasi left continuity of $Y^{(\infty)}$,
     \[
     \lim_{k\to\infty}\E f(Y^{(\overline{n}_k)}_{t_{n_k}})=\E f(Y^{(\infty)}_t),
     \]
     for every bounded continuous function $f$ (keep in mind that the functional tightness of $(M^{(n)})$
follows from Theorem 1.13 in~\cite{JacodShi} which in turn relies on the so-called Aldous criterion; any weak
limiting process of such a sequence in the Skorokhod sense is then quasi-left continuous and so is $Y$ since $B$
is pathwise continuous). Hence 
     $\displaystyle \lim_{k\to\infty}\E f(Y_{n_k})=\E f(Y^{(\infty)}_t)$.
     Observe that the law of $Y^{(\infty)}_0$ is a weak limit of the 
     sequence $Y_n$, so that 
     $\P(Y^{(\infty)}_0\in K^c)<\e$. Now we have
     \[
     \E f(Y_{n_k})-\int f d\nu=\E f(Y_{n_k})-\E f(Y^{(\infty)}_t)
                         +\E f(Y^{(\infty)}_t)-\int f d\nu,
     \]
     so that, if $\mu$ denotes the law of $Y^{(\infty)}_0$,
     \begin{eqnarray*}
     \limsup_{k\to\infty}\left|\E f(Y_{n_k})-\int f d\nu\right|
       &\leq &\left|\E f(Y^{(\infty)}_t)-\int f d\nu\right|\\
       &=&\left|\int \E_y(f(Y_t))d\mu(y)-\int f 
       d\nu\right|\\
       &\leq &\e+2||f||_\infty\mu(K^c)\\
       &\leq &\e(1+2||f||_\infty).
     \end{eqnarray*}
     It follows that any weak limit of the sequence $(Y_n)_{n\in\N}$
     is equal to $\nu$, which completes the proof of Theorem~\ref{TCL}.
     \hfill$\cqfd$
     
     \bigskip
     
     For the proof of Theorem~\ref{stat-dist}, we first observe that
     the generator $L$ depends in an affine way on the state 
     variable $y$. This {\em affine} structure suggests that
     the Laplace transform $\E_y e^{-pY_t}$ has the form
     $e^{\vfi_p(t)+y\psi_p(t)}$, for some functions $\vfi_p$
     and $\psi_p$. Affine models have been recently extensively
     studied in connection with interest rate modelling
     (see for instance~\cite{Duffie} or~\cite{Duffie2}). The following 
     proposition gives a precise description of the Laplace transform.
     \begin{Pro}\label{Laplace}
         Let $(Y_t)_{t\geq 0}$ be the Markov process with generator $L$
         on $[0,+\infty)$.
         We have, for $p>0$, $y\in [0,+\infty)$,
        \begin{equation}\label{YPhiPsi}
        \E_y \,e^{-pY_t}=\exp\left(\vfi_p(t)+y\psi_p(t)\right),
        \end{equation}
        where $\psi_p$ is the unique solution, on $[0,+\infty)$ of the differential 
        equation
        \[
        \psi'=p_{\!_B}\frac{e^{g\psi}-1}{g}-p_{\!_A}\psi, \mbox{ with } 
        \psi(0)=-p,
        \]
        and 
        \[
        \vfi_p(t)=(1-p_{\!_A})\int_0^t\psi_p(s)ds.
        \]
     \end{Pro}
     Before proving the Proposition, we study the involved ordinary 
     differential equation.
      \begin{Lem}\label{ODE}
	Given $\psi_0\in(-\infty,0]$, the ordinary differential
	equation
	\begin{equation}\label{diff-eq}
	\psi'=p_{\!_B}\frac{e^{g\psi}-1}{g}-p_{\!_A}\psi
	\end{equation}
	has a unique equation on $[0,+\infty)$ satisfying
	the initial condition $\psi(0)=\psi_0$. 
	Moreover, we have
	\[
	\forall t\geq 0,\quad \psi(0)\leq \psi(t)e^{\pi t}\leq 0.
	\]
     \end{Lem}
     \noindent{\sc Proof:} Existence and uniqueness of a local
     solution follows from the Cauchy-Lipschitz theorem. 
     In order to prove non-explosion, observe that if 
     $\psi$ solves~(\ref{diff-eq}), we have, using the inequality
     $(e^{g\psi}-1)/g\geq \psi$,
     \[
     \psi'+\pi\psi\geq 0.
     \]
     Therefore, the function $t\mapsto \psi(t)e^{\pi t}$ is 
     non-decreasing, so that $\psi(0)\leq \psi(t)e^{\pi t}$. 
     Since 0 is an equilibrium of the equation, we have
     $\psi(t)\leq 0$ if $\psi(0)\leq 0$, and the inequality is strict 
     unless $\psi(0)=0$. Hence $\psi(0)\leq \psi(t)e^{\pi t}\leq 0$
     and the lemma follows easily.
     \hfill$\cqfd$
     
     \bigskip
     
     {\noindent \sc Proof of Proposition~\ref{Laplace}: }
     Let $u_p(t,y)=\exp(\vfi_p(t)+y\psi_p(t))$, where 
     $\psi_p$ and $\vfi_p$ are defined as in the statement of the 
     Proposition. The existence of $\psi_p$ follows from
     Lemma~\ref{ODE}. An easy computation shows that
     $\frac{\partial u_p}{\partial t}-Lu_p=0$ on 
     $[0,+\infty)\times[0,+\infty)$, so that, for $T>0$,
     the process
     $\left(u_p(T-t,Y_t)\right)_{0\leq t\leq T}$ is a martingale,
     and $\E\, u_p(T,Y_0)=\E\, u_p(0,Y_T)$, and the Proposition follows
     easily.\hfill$\cqfd$
     
     \bigskip
     
     {\noindent \sc Proof of Theorem~\ref{stat-dist}: }
     
     \noindent$\bullet$
     Uniqueness of the invariant distribution. We deduce from Lemma~\ref{ODE} that, with the notation
     of Proposition~\ref{Laplace},  $|\psi_p(t)|\leq e^{-\pi t}$
     and $\displaystyle \lim_{t\to 
     \infty}\vfi_p(t)=(1-p_{\!_A})\int_0^{+\infty}\psi_p(s)ds$. Therefore
     \[
     \lim_{t\to\infty}\E_y(e^{-pY_t})=\exp\left((1-p_{\!_A})\int_0^\infty 
            \psi_p(s)ds \right),
     \]
     and the convergence is uniform on compact sets. This implies
     the uniqueness of the stationary distribution as well 
     as~(\ref{unif}). We also have the Laplace transform of $\nu$:
     \[
     \int_{\R^+} e^{-py}\nu(dy)=\exp\left((1-p_{\!_A})\int_0^\infty 
            \psi_p(s)ds \right).
     \]
     Note that, since $\psi_p\leq 0$ and 
     $\psi_p'=p_{\!_B}\frac{e^{g\psi_p}-1}{g}-p_{\!_A}\psi_p$, we have 
     $\psi_p'+p_{\!_A}\psi_p\leq 0$. Therefore, $\psi_p(t)\leq 
     -pe^{-p_{\!_A}t}$,
     and
     \[
    \forall p\geq 0,\quad \int e^{-py}\nu(dy)\geq \exp(-p(1-p_{\!_A})/p_{\!_A})=\exp(-pr_{\!_A}).
     \]
     This yields $\nu([0,r_{\!_A}))=0$.
     
     \ss
\noindent$\bullet$ Further properties of the invariant distribution $\nu$. 
     The stationary distribution satisfies $\int Lf d\nu =0$ for any 
     continuously differentiable function $f$ with compact support
     in $[0,+\infty)$.
     This reads
     \begin{equation}\label{nuLf}
     \forall f\in C^1_{K}, \quad \int \left(r y\frac{f(y+g)-f(y)}{g}
          +(r_{\!_A}-y)f'(y)\right)\!\nu(dy)=0,
     \end{equation}
     where $r=p_{\!_B}/p_{\!_A}$ and $r_{\!_A}=(1-p_{\!_A})/p_{\!_A}$.
     
     We first show that $\nu(\{r_{\!_A}\})=0$. Let $\vfi$ be a non-negative
     continuously differentiable function satisfying $\vfi=1$ in a 
     neighbourhood of the origin and $\vfi=0$ outside the interval
     $[-1,1]$. For $n\geq 1$ let
     \[
     f_n(y)=\vfi(n(y-r_{\!_A})), \quad y\in\R.
     \]
     We have $f_n(y)=0$ if $|y-r_{\!_A}|\geq 1/n$. In particular, the 
     support of $f_n$ lies in $[0,+\infty)$, for $n$ large enough. Applying (\ref{nuLf})
     with $f=f_n$, we get
     \[
     \int \left(r y\frac{f_n(y+g)-f_n(y)}{g}
          +(r_{\!_A}-y)n\vfi'(n(y-r_{\!_A}))\right)\!\nu(dy)=0.
     \]
     Observe that $\displaystyle \lim_{n\to\infty}f_n=\ind_{\{r_{\!_A}\}}$ so that
     \[
     \lim_{n\to\infty}\int  y(f_n(y+g)-f_n(y))\nu(dy)=
     (r_{\!_A}-g)\nu(\{r_{\!_A}-g\})-r_{\!_A}\nu(\{r_{\!_A}\})=-r_{\!_A}\nu(\{r_{\!_A}\}),
     \]
     where we have used $\nu(-\infty,r_{\!_A})=0$. On the other hand,
     we have $|(r_{\!_A}-y)n\,\vfi'(n(y-r_{\!_A}))|\leq \sup_{u\in\R}(u\vfi'(u))$,
     and $\displaystyle \lim_{n\to\infty}(n\vfi'(n(y-r_{\!_A})))=0$, so that, by dominated convergence,
     \[
     \lim_{n\to\infty}\int (r_A-y)n\,\vfi'(n(y-r_{\!_A}))\nu(dy)=0.
     \]
     Hence $\nu(\{r_{\!_A}\})=0$.
     
     We now study the measure $\nu$ on the open interval 
     $(r_{\!_A},+\infty)$. Denote by ${\cal D}$ the set of all infinitely
     differentiable functions with compact support in $(r_{\!_A},+\infty)$.
     We deduce from (\ref{nuLf}) that, for $f\in{\cal D}$,
     \begin{equation}\label{nuLf1}
     \frac{r}{g}\int \nu(dy)y f(y+g)-\frac{r}{g}\int \nu(dy)y f(y)
     +\int \nu(dy)(r_{\!_A}-y)f'(y)=0.
     \end{equation}
     Denote by $\nu_g$ the measure defined by 
     $ \displaystyle \int\nu_g(dy)f(y)=\int \nu(dy)f(y+g)$. We deduce 
     from~(\ref{nuLf1}) that
     $\nu$ satisfies the following equation in the sense of 
     distributions:
     \[
     (y-r_{\!_A})\nu'+(1-(r/g)y)\nu=-\frac{r}{g}(y-g)\nu_g,
     \]
     or
     \begin{equation}\label{equanu}
     \nu'+\frac{1-(r/g)y}{y-r_{\!_A}}\nu=-\frac{r}{g}\frac{y-g}{y-r_{\!_A}}\nu_g.
     \end{equation}
     Denote by $F$ the function defined by
     \begin{equation}\label{F}
     F(y)=e^{ry/g}(y-r_{\!_A})^{d-1},\quad y>r_{\!_A},
       \end{equation}
     where $d=r\,r_{\!_A}/g$. We have
     \[
     F'(y)=-\frac{1-(r/g)y}{y-r_{\!_A}}F(y),
     \]
     so that the equation satisfied by $\nu$ reads
     \begin{equation}\label{nu}
         \left(\frac{1}{F}\nu\right)'=\frac{G}{F}\nu_g,
     \end{equation}
     where the function $G$ is defined by $G(y)=-\frac{r}{g}\frac{y-g}{y-r_{\!_A}}$.
     
     On the set $(r_{\!_A},r_{\!_A}+g)$, the measure $\nu_g$ vanishes, so that
     $\nu=\lambda_0 F$ for some non negative constant $\lambda_0$.
     At this point, we know that the restriction of the measure $\nu$
     to the set $(0, r_{\!_A}+g)$ has a density  which vanishes on
     $(0,r_{\!_A})$ and is given by $\lambda_0 F$ on $(r_{\!_A},r_{\!_A}+g)$.
     
     \medskip
     
     We will prove by induction that the distribution  $\nu$  coincides 
     with a continuous function on $(r_{\!_A},r_{\!_A}+ng)$, which is infinitely
     differentiable on $(r_{\!_A}+(n-1)g, r_{\!_A}+ng)$.
     The claim has been proved for $n=1$. Assume that it is true for
     $n$. On the set $(r_{\!_A},r_{\!_A}+(n+1)g)$, the distributional derivative
     of $(1/F)\nu$ coincides with the function $y\mapsto 
     (G(y)/F(y))\nu(y-g)$, which is locally integrable  on 
     $(r_{\!_A},r_{\!_A}+ng+g)$, continuous on $(r_{\!_A}+g, r_{\!_A}+ng+g)$,
     and infinitely differentiable on $(r_{\!_A}+ng,r_{\!_A}+ng+g)$,
     due to the induction hypothesis (there may be a discontinuity at
     $r_{\!_A}+g$ if $d<1$). It follows that $(1/F)\nu$ is a continuous 
     (resp. infinitely differentiable)
     function, and so is $\nu$ on  $(r_{\!_A},r_{\!_A}+(n+1)g)$ (resp. 
     $(r_{\!_A}+ng,r_{\!_A}+ng+g)$).
     We have proved that $\nu$ has a continuous density on 
     $(r_{\!_A},+\infty)$, which is infinitely differentiable on the open
     set $\bigcup_{n=1}^\infty (r_{\!_A}+(n-1)g, r_{\!_A}+ng)$.
     
     \medskip
     
     Finally, we  prove that the density of $\nu$ is positive on 
     $(r_{\!_A},+\infty)$. Note that $G(y)<0$ if $y>g$ and that the density
     vanishes at $y-g$ if $y<g$. Therefore 
     $\left(\frac{1}{F}\nu\right)'\leq 0$, so that the function
     $y\mapsto \nu(y)/F(y)$ is nondecreasing. 
     It follows that $\lambda_0$ cannot be zero (otherwise $\nu$
     would be identically zero). Hence $\nu(y)>0$ for 
     $y\in(r_{\!_A},r_{\!_A}+g)$. Now, if $\nu(y)>0$ for $y\in(r_{\!_A}+ng-g,r_{\!_A}+ng)$,
     the function $\nu/F$ is strictly decreasing on $(r_{\!_A}+ng,r_{\!_A}+ng+g)$
     and, therefore, cannot vanish. So, by induction, the density is positive on 
     $(r_{\!_A},+\infty)$. This completes the proof of 
     Theorem~\ref{stat-dist}.\hfill$\cqfd$

\bigskip
\noindent {\bf Additional remarks.} $\bullet$ The proof of Theorem~\ref{stat-dist}  provides a bit more information on the
invariant distribution   $\nu$. Let $g>0$ and let  $\phi_g$ denote   its continuous density on
$(r_{\!_A},+\infty)$:  the function $\phi_g$ is ${\cal C}^\infty$ on
$[r_{\!_A},+\infty)\setminus(r_{\!_A}+g\,\N)$ and it follows from~(\ref{F}) and the definitions of $r$ and
$r_{\!_A}$ (and $d=rr_{_A}/g$, see the proof of theorem~\ref{stat-dist}) that
\[
\phi_g(r_{\!_A})=+\infty \;\mbox{ if }\; g>g^*,\quad \phi_g(r_{\!_A})\!\in(0,+\infty)
\; \mbox{ if }\; g=g^*\;\mbox{ and }\; \phi_g(r_{\!_A})=0 
\; \mbox{ if }\; g<g^*
\] 
 where $g^*= \frac{p_{\!_B}(1-p_{\!_A})}{p_{\!_A}^2}\!\in(0,
\frac{1-p_{\!_A}}{p_{\!_A}})$. As concerns the regularity of the density
$\phi_g$ at points
$y \in r_{_A}+g\,\N$, one easily derives from Equation~(\ref{equanu}) that
for every $m,\, k\in\N$, 

\smallskip
-- $\phi_g$ is $C^{m+k}$ at $r_{_A}+kg$ as soon as
$g<\frac{g^*}{m+1}$, 

-- the $(m+k)^{th}$ derivative $\phi_g^{(m+k)}$ is only  right and
left continuous at
$r_{_A}+kg$ if $g=\frac{g^*}{m+1}$.    

\medskip
\noindent $\bullet$ One can  characterize the finite positive exponential moments of $\nu$ by slightly exten\-ding
the proof of  Proposition~\ref{Laplace} (Laplace transform). For every $y>1$, let  $\theta(y)$ denote the unique
(strictly) positive solution of the equation 
\[
\frac{e^\theta-1}{\theta}=y.
\] 
Note  that $\log y< \theta(y) <2(y-1)$ and that $\displaystyle \lim_{y\to 1} \frac{\theta(y)}{2(y-1)}=1$ and
$\displaystyle \lim_{y\to \infty} \frac{\theta(y)}{\log y}=1$. The result is as follows 
\begin{equation}\label{Laplaceposit}
 \int e^{py}\nu(dy) <+\infty \quad \mbox{ if and only if }\quad p < p^*_g:=g\, \theta(p_{\!_A}/p_{\!_B}).
\end{equation}

With the notations of Proposition~\ref{Laplace}, it follows from Fatou's  Lemma that 
\begin{equation}\label{Fatounu}
\forall\, p>0,\qquad \int e^{py}\nu(dy)\le \liminf_{t\to \infty} \E_y (e^{pY_t}).
\end{equation}
We know that
\[
\E_y (e^{pY_t}) = e^{\tilde{\varphi}_p(t) +y\tilde{\psi}_p(t)}  
\]    
with   $\tilde{\varphi}_p(t) =
(1-p_{\!_A})\int_0^t\tilde{\psi}_p(s)ds$ and $\tilde{\psi}_p$ is solution on the non-negative real line (if any) of 
\[
\psi'(t) =G(\psi(t)),\quad\psi(0)=p\quad  \mbox{ with }\quad G(u) = -p_{\!_A}u+\frac{p_{\!_B}}{g}(e^{gu}-1).
\]
The function $G$ is convex on $\R_+$  and satisfies $G(0)=G(p^*_g)=0$, $G((0,p^*_g))\subset (-\infty,0)$. 

Let
$p\!\in(0,p^*_g)$. The convexity of $G$ implies 
\[
\forall\, u\in [0,p], \quad \frac{G(u)}{u} \le \frac{G(p)}{p}<0.
\] 
It follows that $\tilde{\psi}_{p}$ does exist on $\R_+$ and satisfies $\displaystyle 0\le \tilde{\psi}_p(t) \le
pe^{\frac{G(p)t}{p}}$ (hence it goes to $0$ when $t$ goes to infinity). One derives that 
\[
\lim_{t\to +\infty} \tilde{\varphi}_p(t)=(1-p_{\!_A})\int_0^{+\infty} \tilde{\psi}_p(t)dt  \le -(1-p_{\!_A})\frac{p^2}{G(p)}.
\]
Combining this with~(\ref{Fatounu}) yields  
\[
\int e^{py}\nu(dy)\le e^{-(1-p_{\!_A})\frac{p^2}{G(p)}}<+\infty.
\]
On the other hand if $p=p^*_g$, $\tilde{\psi}_p(t) = p^*_g$ and $\tilde{\varphi}_p(t) = (1-p_{\!_A})p^*_gt$. Consequently
\[
\forall\, t\ge 0,\qquad \int e^{p^*_gy}\nu(dy) = \int \E_y(e^{p^*_gY_t})\nu(dy)= e^{(1-p_{\!_A})p^*_gt} \int
e^{p^*_gy}\nu(dy).
\]  
Now the right hand side of this equality goes to $\infty$ as $t$ goes to infinity since $(1-p_{\!_A})p^*_g\!>\!0$ which
shows that $\displaystyle \int\! e^{p^*_gy}\nu(dy)\!=\!+\infty$ (since it cannot be $0$).

\medskip
\noindent $\bullet$ One has, in accordance with  the convergence rate result obtained for  $\rho_n=o(\g_n)$, that 
\[
\int y\,\nu(dy) = \frac{1-p_{\!_A}}{\pi}.
\]
To prove this claim, one first notes, using  the definition~(\ref{generateur}) of the generator $L, $ that
$L(Id)(y) = 1-p_{\!_A}-\pi\,y$. Hence the above claim will follow from
$\displaystyle \int L(Id)(y)\nu (dy)=0$.  Let $\varphi:\R_+\to \R_+$ denote a
continuously differentiable function such that $\varphi(y)=y$ if $y\!\in [0,1]$,
$\varphi(y)=0$ if $y\ge 2$ and $\varphi'$ is bounded on $\R_+$. Set $\varphi_n(y)
= n\varphi(y/n)$, $n\ge 1$. One checks that $L(\varphi_n)\to L(Id)$ as $n$ goes
to infinity  and $|L(\varphi_n)(y)|\le ay+b$ for some positive real constants
$a,\,b$. One derives by the dominated convergence theorem that 
\[
\int L(Id)(y)\nu(dy) =\lim_n\int L(\varphi_n)(y)\nu(dy) = 0 
\]
where we used that the function $\varphi_n$ has compact support on $[0,+\infty)$.  One shows similarly that 
 $\displaystyle \int
\!L(u\mapsto u^2)(y)\nu(dy)=0$ to derive that
\[
\int \left(y-\frac{1-p_{\!_A}}{\pi}\right)^2\nu(dy) =
g\,\frac{p_{\!_B}(1-p_{\!_A})}{2\pi^2}. 
\]
Note that, as one could expect,  this variance goes to $0$ as $g\to 0$. As a conclusion, we present in figure~1
three  examples of shape for $\phi_g$. They were obtained from an exact simulation of the Markov process
$(Y_t)_{t\ge 0}$ (associated to the     generator 
$L$) at its jump times: we approximated the p.d.f.  by a histogram method using 
Birkhoff's  ergodic Theorem.

\bs
\centerline{Figures should be here}

\begin{figure}
\begin{center}
%\begin{tabular}{ccc}  
%$$\includegraphics[height =4 cm,width=4cm]{Fig1(erg).eps}$$&  
%$$\includegraphics[ height =4 cm, width=4cm]{Fig2(erg).eps}$$& 
%$$\includegraphics[height =4 cm,width=4cm]{Fig3(erg).eps}$$  
% \end{tabular}
\caption{\em   Graphs of the p.d.f $\phi_g$, $p_{_{\!A}}= 2/5$, $g=1$; the vertical dotted line shows the mean
$\frac{1-p_{_{\!A}}}{\pi}$ of
$\nu$. {\rm Left:}
$p_{_{\!B}}=1/3$ ($g^*>g=1$). {\rm Center:} 
$p_{_{\!B}}=4/15$ ($ g^*=g=1$). {\rm Right:} $p_{_{\!B}}=1/6$ ($g^*<g=1$).}
%\label{fig:AsCall_B}
\end{center}
\end{figure}

\bs\noindent {\bf A final remark about the case $\pi=0$ and
$\g_n=g\,\rho_n$.} In that setting (see~Remark~\ref{Rem1}) the asymptotics
of the algorithm cannot be elucidated by using the $ODE$ approach since
it holds in a weak sense. Setting $Y_n = 1-2X_n$ one checks that $Y_n\!\in[-1,1]$ and 
\[
Y_{n+1} =Y_n (1-2g\rho^2_{n+1}(1-p_{_{\!A}}))-2g\,\rho_{n+1}\Delta M_{n+1}
\]
and that $\E((\Delta M_{n+1})^2\,|{\cal F}_{n+1})= \frac{p_{_{\!A}}}{4}(1-Y^2_n)+O(\rho^2_{n+1})$.
Then, a similar approach as that developed in this section (but significantly less technical since $(Y_n)$
is bounded by $1$) shows that $Y_n$ converges in distribution to the
invariant distribution $\mu$ of the Brownian diffusion with generator
${\cal L}f(y) = -2g(1-p_{_{\!A}})yf'(y)+\frac 12g^2p_{_{\!A}}(1-y^2) f''(y)$. In that case, it is
well-known that $\mu$ has a density function for which a closed form is
available (see~\cite{KATA}), namely
\[
\mu(dy)= m(y)dy\quad \mbox{ with }\quad m(y) = C_{g,r_{_A}}
(1-y^2)^{\frac{2r_{A}}{g}-1}\mbox{\bf 1}_{(-1,1)}(y).
\]
Note that when  
$g=2r_{_{\!A}}=2(1/p_{_{\!A}}-1)>0$, $\mu$ is but the uniform distribution over $[-1,1]$.

\end{document}